\RequirePackage{fix-cm}
\documentclass[smallextended]{svjour3}       % onecolumn (second format) % smallextended
\smartqed  % flush right qed marks, e.g. at end of proof
\usepackage[utf8]{inputenc} %unicode support
\usepackage{lineno,hyperref}
\usepackage{amsmath, amssymb}
\modulolinenumbers[5]
\usepackage{lscape}
\usepackage{graphicx}
\usepackage{booktabs}
\usepackage{epstopdf}
\usepackage{float}
\usepackage{color} %{\color{blue}   }
\usepackage{cases}
\usepackage{longtable}
\usepackage{multirow}
\usepackage[figuresright]{rotating} % 制作横向表格
\usepackage{algorithmic}
\usepackage[ruled,linesnumbered]{algorithm2e}
\usepackage[figuresright]{rotating}
\usepackage{marvosym} %通讯作者小信封

\newtheorem{assumption}{Assumption}

\newcommand{\ba}{\begin{array}}
\newcommand{\ea}{\end{array}}
\newcommand{\be}{\begin{equation}}
\newcommand{\ee}{\end{equation}}

\newlength\savewidth

\begin{document}

\title{Monotone Splitting SQP Algorithms for Two-block Nonconvex Optimization  Problems with General Linear Constraints and Applications
\thanks{ This research was supported by the National Natural Science Foundation of China (12171106 and 12261008), and the Guangxi Natural Science Foundation (2020GXNSFDA238017 and 2018GXNSFFA281007) and the Xiangsihu Young Scholars Innovative Research Team of Guangxi Minzu University
(2022GXUNXSHQN04).}
}

\titlerunning{Monotone Splitting SQP Algorithms}        % if too long for running head

\author{Jinbao Jian$^{1}$          \and
        Guodong Ma$^{1}$          \and
        Xiao Xu$^{2}$          \and
        Daolan Han$^{1}$         % \and
%        Lin Yang$^{1}$  %etc.
}

\authorrunning{Jinbao Jian et al.} % if too long for running head

%\institute{F. Author \at
%              first address \\
%              Tel.: +123-45-678910\\
%              Fax: +123-45-678910\\
%              \email{fauthor@example.com}           %  \\
%%             \emph{Present address:} of F. Author  %  if needed

\institute{
           Jinbao Jian  \at
%              Fax: +123-45-678910\\
              \email{jianjb@gxu.edu.cn}           %  \\
           \and
            \Letter Guodong Ma   \at
%              Fax: +123-45-678910\\
              \email{mgd2006@163.com}           %  \\
           \and
           Xiao Xu  \at
%              Fax: +123-45-678910\\
              \email{1209077536@qq.com}           %  \\
           \and
          Daolan Han \at
%              Fax: +123-45-678910\\
              \email{handaolan@126.com}                 %  \\
         %  \and
%           Lin Yang    \at
%              %Tel.: +123-45-678910\\
%%              Fax: +123-45-678910\\
%              \email{yanglinlin1125@163.com}          %  \\
           \and
              1. College of Mathematics and Physics, Center for Applied Mathematics of Guangxi, Guangxi Key Laboratory of Hybrid Computation and IC Design Analysis, Guangxi Minzu University, Nanning, Guangxi, 530006, China \\
              2. College of Mathematics and Information Science, Guangxi University, Nanning, Guangxi, 530004, China
}

% AMS subject classifications: 90C26, 90C30, 65K05

\maketitle

\begin{abstract}
In this work, based on the ideas of alternating direction method with multipliers (ADMM) and sequential quadratic programming (SQP), as well as Armijo line search technology, monotone splitting SQP algorithms for two-block nonconvex optimization problems with linear  equality, inequality and box constraints are discussed. Firstly, the discussed problem is transformed into an optimization problem with only linear equality and box constraints by introducing slack variables. Secondly, we use the idea of ADMM to decompose the quadratic programming (QP)~subproblem. Especially, the QP subproblem corresponding to the introducing slack variable is simple, and it has an explicit optimal solution without increasing computational cost. Thirdly, the search direction is generated by the optimal solutions of the subproblems, and the new iteration point is yielded by Armijo~line search with augmented Lagrange function. And the global convergence of the algorithm is analyzed under weaker assumptions. In addition, box constraints are extended to general nonempty closed convex sets, moreover, the global convergence of the corresponding algorithm is also proved. Finally, some
preliminary numerical
experiments and applications in the mid-to-large-scale economic dispatch problems for power systems are reported, and these show that our proposed algorithm is promising.

\end{abstract}

\section{Introduction}
\label{intro}

Let us start with the canonical two-block separable  optimization problem with linear
equality constraints,
\be\label{1.10x}
\{\min\ f(x)+\theta(y)|\ Ax+By=b,\ x\in \mathcal{X},\ y\in \mathcal{X}\},
\ee
where function\ $f:\mathfrak{R}^{n_{1}}\rightarrow \mathfrak {R}$  is proper and lower semi-continuous,~$\theta:\mathfrak {R}^{n_{2}}\rightarrow \mathfrak {R}$ is a continuous differentiable function, $A\in \mathfrak {R}^{m_{1}\times n_{1}},~B\in \mathfrak {R}^{m_{1}\times n_{2}},$ $b\in \mathbb{R}^{m_1}$ are the given matrices and vector.

Many  problems  can be expressed in the form of problem (\ref{1.10x}),
such as data mining, signal and image processing, electric power
systems, etc, e.g., \cite{WLFK,BPCP,XCXZ,ZYJ}. Because of its separable structure, problem (\ref{1.10x}) can be efficiently solved by the Douglas-Rachford (DR) splitting method and the Peaceman-Rachford (PR) splitting method (PRSM).

Although the augmented Lagrangian method (ALM) can be applied to solve problem(\ref{1.10x}), it does
not take full advantage of the separable structure of (\ref{1.10x}).
As a splitting version of ALM,
the standard alternating direction method of multipliers (ADMM, \cite{GM,GB}), which can be viewed as an application of DR splitting method solving the augmented Lagrangian function (ALF), exploits the
separable structure of the objective function and performs the following iterations:

%[11] R. Glowinski and A. Marrocco, Approximation pare′le′ments finis d’rdre un et re′solution, par
%pe′nalisation-dualite′ d’une classe de proble`mes de Dirichlet non line′aires, Rev. Fr. Autom.
%Inform. Rech. Ope′r. Anal. Nume′r., 2 (1975), 41-76.

%The alternating direction method with multipliers has been widely researched since 1970s,. For two block convex problems, its convergence has been well analyzed.
%The author introduced ADMM in \cite{DB1976,RA1975}, and established its convergent result for strongly convex functions. The scholars extended associated result to general convex functions(see \cite{J1989,JD1992}). We know that
%ADMM can converge at a sublinear convergence rate of $\mathcal{O}(1/n)$\cite{BX2012,RB2013}, and $\mathcal{O}(1/n^{2})$ for the accelerated version \cite{TBS2012}. There are other variations of ADMM for convex objective functions, i.e.\cite{FZH2014,FWZ2015,GT2015}. However, the convergence analysis of ADMM has always been an open question for nonconvex problems. Fortunately, scholars have done some research on nonconvex ADMM, i.e.,\cite{FZH2014,GT2015,KDT2016,MZM2014}.
%The procedure of classical ADMM for (\ref{1.10x}) is as follows:
\be\label{1.111y}
\left\{\begin{array}{lll}\nonumber
x^{k+1}\in \arg\min\{\mathcal{L}_\beta(x,y^{k},\lambda^{k})|\ x\in \mathcal{X}\},\\
y^{k+1}\in \arg\min\{\mathcal{L}_\beta(x^{k+1},y,\lambda^{k})|\ y\in \mathcal{Y}\},\\
\lambda^{k+1}=\lambda^{k}-s\beta(Ax^{k+1}+By^{k+1}-b),
\end{array}
    \right.
    \ee
where $s\in(0,\frac{1+\sqrt{5}}2)$ is the stepsize for updating the dual variable $\lambda$, $\beta>0$ is the penalty parameter, $\mathcal{L_\beta(\cdot)}$ is the ALF of (\ref{1.10x}) and defined as follows: \be\label{1.2}\mathcal{L}_\beta(x,y,\lambda) = f(x)+\theta(y)-\langle\lambda,Ax+By-b\rangle+\frac{\beta}{2}\|Ax+By-b||^2.
\ee
If the PRSM \cite{PR} is applied to the dual of problem
(\ref{1.10x}), then we obtain a variation of ADMM, whose
iteration scheme is as follows:

\be\label{1.111y}
\left\{\begin{array}{lll}\nonumber
x^{k+1}\in \arg\min\{\mathcal{L}_\beta(x,y^{k},\lambda^{k})|\ x\in \mathcal{X}\},\\
\lambda^{k+\frac12}=\lambda^{k}-\beta(Ax^{k+1}+By^{k}-b),\\
y^{k+1}\in \arg\min\{\mathcal{L}_\beta(x^{k+1},y,\lambda^{k})|\ y\in \mathcal{Y}\},\\
\lambda^{k+1}=\lambda^{k}-\beta(Ax^{k+1}+By^{k+1}-b).
\end{array}
    \right.
    \ee
The PRSM above is also called symmetric ADMM since the Lagrange multipliers are
symmetrically updated twice in each loop. Note that both updates of dual variable in
PRSM use the same constant stepsize 1.
Motivated from the ideas of enlarging the dual stepsize in \cite{HY1}, the following extension of the symmetric ADMM was developed
by He, et al. \cite{HY2}:

\be\label{1.111y}
\left\{\begin{array}{lll}\nonumber
x^{k+1}\in \arg\min\{\mathcal{L}_\beta(x,y^{k},\lambda^{k})|\ x\in \mathcal{X}\},\\
\lambda^{k+\frac12}=\lambda^{k}-r\beta(Ax^{k+1}+By^{k}-b),\\
y^{k+1}\in \arg\min\{\mathcal{L}_\beta(x^{k+1},y,\lambda^{k})|\ y\in \mathcal{Y}\},\\
\lambda^{k+1}=\lambda^{k}-s\beta(Ax^{k+1}+By^{k+1}-b),
\end{array}
    \right.
    \ee
where, for the sake of convergence, the stepsize pair $(r,s)$ is required to belong to the
following region:
$$D_0=\{(r,s)|\ s\in(0,\frac{1+\sqrt{5}}2),r+s>0,r\in(-1,1),\ |r|<1+s-s^2\}. $$

If the linear equality constraints in (\ref{1.10x}) are changed to linear inequality constraints while
all the other settings are remained, we obtain the following model:
\be\label{1.11x}
\min\{ f(x)+\theta(y)|\  Ax+By\geq b,\ x\in \mathcal{X},\ y\in \mathcal{Y}\}.
\ee
The two-block separable % convex
optimization model (\ref{1.11x}) with linear inequality constraints
captures particular applications such as the support vector machine with a linear kernel in \cite{VV}
and its variant in \cite{LMO}. To solve problem (\ref{1.11x}), by introducing an auxiliary variable $z$, it can be reformulated as the
following three-block separable model with linear equality constraints:
\be\label{1.12x}
\min\{ f(x)+\theta(y)|\  Ax+By-z= b,\ \ x\in \mathcal{X},\ y\in \mathcal{Y},\ z\geq0\}.
\ee
Then, a direct extension of the ADMM % (\ref{1.111y}), 对不上号！
 can be
applied to the reformulated model (\ref{1.12x}) resulting in the following iterative scheme:
\be\label{1.112y}
\left\{\begin{array}{lll}\nonumber
x^{k+1}\in\arg\min\{\mathcal{L}_\beta(x,y^{k},z^{k},\lambda^{k})|\ x\in \mathcal{X}\},\\
y^{k+1}\in\arg\min\{\mathcal{L}_\beta(x^{k+1},y,z^{k},\lambda^{k})|\ y\in \mathcal{Y}\},\\
z^{k+1}\in\arg\min\{\mathcal{L}_\beta(x^{k+1},y^{k+1},z,\lambda^{k})|\ z\geq 0 %z\in \in \mathcal{R}_+^{r}
\},\\
\lambda^{k+1}=\lambda^{k}-\beta(Ax^{k+1}+By^{k+1}-z^{k+1}-b),
\end{array}
    \right.
    \ee
more specifically, the ALF $\mathcal{L_\beta(\cdot)}$ of (\ref{1.12x}) is defined as follows: \be\label{1.2x}\mathcal{L}_\beta(x,y,z,\lambda) = f(x)+\theta(y)-\langle\lambda,Ax+By-z-b\rangle+\frac{\beta}{2}\|Ax+By-z-b||^2.
\ee
 According to \cite{CHYY}, however, convergence of the
direct extension of ADMM (\ref{1.112y}) is not guaranteed unless more restrictive conditions on the
objective functions, coefficient matrices, as well as the penalty parameter, are additionally posed.
Alternatively, the iterative scheme (\ref{1.112y}) should be revised appropriately to render the convergence. For
example, the scheme  (\ref{1.112y}) should be further corrected by those correction steps studied in
\cite{HTY,HY3}. Recently, He, Xu and Yuan \cite{HXY} proposed a unified framework of algorithmic design and a roadmap for convergence
analysis on the extensions of ADMM for separable convex optimization problems with linear equality or inequality constraints.

On the other hand, it is known that the sequential quadratic
programming (SQP) method is a very
important technique for designing efficient algorithms for smooth constrained optimization problems \cite{a18xx,a18,Solodov2009,jianjinbao_2010,a18xxx,z3}. In recent years, in order to better develop SQP algorithms for separable nonconvex optimization problems,  Jian et al. introduced
the idea of a splitting method for solving quadratic
programming (QP) subproblems, namely, splitting
the large-scale QP into two or more small-scale QP subproblems. As a
result, a class of splitting SQP algorithms are
proposed; see, e.g., \cite{a20,a21,a19}.

The monotone splitting SQO
algorithm \cite{a20}, problem $\min\{f(x)+\theta(y) |\ Ax + By = b, x \in [l, u], y \in [p, q]\}$ is discussed, and the main
features of this algorithm are as follows: First, the primal search direction associated with the
primal variable is yielded by solving two independent small-scale QP
subproblems in parallel, a deflection
of the steepest descent direction of the ALF (\ref{1.2}) for the dual variable is chosen as
the search direction for the dual variable. Second,  the primal-dual variables
are considered as a whole, ALF is the merit function, and a new primal-dual iterative point is generated by Armijo line search.

For the two-block nonconvex optimization problem  $\min\{f(x)+\theta(y) |\ h(x) + g(y) = 0, x \in [l, u], y \in [p, q]\}$ , Jian et al. \cite{a21} proposed a quadratically constrained quadratic optimization (QCQO)-based splitting SQQ algorithm.
The basic ideas in \cite{a21} can be summarized as follows:
First, a QCQO subproblem for the discussed problem at the current iteration is considered.
Second, with the help of ALF dealing with equality constraints and splitting techniques, the QCQO problem is split into two small-scale subproblems with nonquadratic objective and affine-box constraints.
Third, a new primal iterative point is generated by the Armijo line search, ALF, as a merit function along the obtained improved direction. The dual multiplier variable is yielded in a similar pattern to algorithm \cite{a20}.

%The Peaceman-Rachford (PR) splitting  SQP algorithm \cite{a19} discusses problem $\min\{f(x)+\theta(y) |\ Ax + By = b, l\leq Cx \leq u, \ p\leq Dy\leq q\}$,
%the basic ideas in \cite{a19} can be summarized
%as follows: First,  based on the idea of the PR splitting algorithm, the augmented Lagrangian
%quadratic programming  in the classical SQP method is decomposed into two small-scale QPs, then the improved search directions are obtained by solving two small-scale QPs; Second, by taking the ALF as the merit function, Armijo line searches are performed along the two improved directions such that
%two step-lengths are yielded, this line search technique not only guarantees the global and strong convergence
%as well as reasonable iteration complexity of the method, but also overcomes the Maratos effect; Third, a new symmetric update technique for the multiplier associated with the equality constraints is
%presented.

In this work, motivated by the idea of ADMM and splitting SQP algorithm for separable optimization problem in \cite{a20}, we propose monotone splitting SQP algorithm for two-block
nonconvex optimization problems with linear equality, inequality and box constraints. And our work possesses the following features:

(i) The discussed problem is transformed into an optimization problem with only linear equality and box constraints by introducing slack variables;

(ii) Use the idea of ADMM and splitting SQP algorithm to decompose the QP~subproblem. Especially, the QP subproblem corresponding to the introducing slack variable is simple, and it has an explicit optimal solution without increasing computational cost;

(iii) The search direction is yielded by the optimal solutions of the subproblems, and the new iteration point is generated by Armijo~line search with the ALF;

(iv) And the global convergence of our proposed algorithm is analyzed. In addition, box constraints are extended to general nonempty closed convex sets, moreover, the global convergence of the corresponding algorithm is also proved.

The paper is organized as follows. The next section describes the motivation and algorithm. Sections 3 and 4 discusses the convergence % analysis
 and extension of the method, respectively. Section 5 contains applications
in the electric power systems. Finally, conclusions
are given in Section 6.

\textbf{Notation.} For any column vectors $x, y, \ldots, \lambda,$ and matrix $C$, throughout this paper, we denote that $(x,y,\lambda,\ldots) :=(x^\top,y^\top,\lambda^\top,\ldots)^\top$, and $||x||^2_C :=x^\top Cx$. $p\bot q$ represents $p^\top q=0$. $||x||$ and $||C||$ represent the $\ell_2$ norms of vector $x$ and matrix $C$, respectively. $C\succ 0$ represents that the matrix $C$ is symmetric positive definite, and $C\succ D$ represents that the matrix $C-D$ is positive definite.

\section{Motivation and algorithm}
\label{sec:1}

In this work, we consider the general linear constrained two-block nonconvex optimization problem with the following form:
\be\label{1.1}
\ba{ll} \min &f(x)+\theta(y)\\
{\rm s.t.}&Ax+By=b,\\
&r\leq Cx+Dy\leq s,\\
&x\in \mathcal{X}:=[l,u], y\in \mathcal{Y}:=[p,q],
\ea
\ee
where~$f:\mathfrak{R}^{n_{1}}\rightarrow \mathfrak {R}$ and~$\theta:\mathfrak {R}^{n_{2}}\rightarrow \mathfrak {R}$ are both smooth but not necessarily convex, and $A\in \mathfrak {R}^{m_{1}\times n_{1}},~B\in \mathfrak {R}^{m_{1}\times n_{2}},$ $C\in \mathfrak {R}^{m_{2}\times n_{1}},D\in \mathfrak {R}^{m_{2}\times n_{2}},$ $b\in \mathbb{R}^{m_1}$ are the given matrices and vector. The box constraints~$l\in \mathfrak {R}^{n_{1}}\cup\{-\infty\}^{n_1},u\in \mathfrak {R}^{n_{1}}\cup\{+\infty\}^{n_1},$ $p\in \mathfrak {R}^{n_{2}}\cup\{-\infty\}^{n_2},q\in \mathfrak {R}^{n_{2}}\cup\{+\infty\}^{n_2},$ without loss of generality, we suppose that $r_i<s_i,\ l_i<u_i$ and $p_i<q_i$.

To present our analysis in a compact way, define matrices
$$\begin {array}{ll} E=\left(\ba{c} A\\C \ea\right), ~F=\left(\ba{c} B\\D \ea\right),~G=\left(\ba{c} 0_{m_{1}\times m_{2}}\\-I_{m_{2}\times m_{2}} \ea\right),~c=\left(\ba{c} b\\0 \ea\right)\end {array}.$$ Then let $Cx+Dy=z\in\mathbb{R}^{m_2}$, and the problem (\ref{1.1}) can be reformulated as follows:
\be\label{1.2}
\ba{ll} \min &f(x)+\theta(y)\\
{\rm s.t.}&Ex+Fy+Gz=c,\\
&l\leq x\leq u,\ p\leq y\leq q,\ r\leq z\leq s.
\ea
\ee
%In order to make the research results more widely applicable, the non-negative restriction of variable~$z$ in problem (\ref{1.2}) can be extended to the generalized box constraint~$[r,s],$ i.e.
%\begin{equation}\label{y0}
%\ba{ll} \min &f(x)+\theta(y)\\
%{\rm s.t.}&Ex+Fy+Gz=c,\\
%&l\leq x\leq u,\ p\leq y\leq q,\ r\leq z\leq s.
%\ea
%\end{equation}
The full Lagrangian function of the problem~(\ref{1.2}) is defined as
\begin{eqnarray}\label{1.3}
&&{L}(x,y,z,\lambda,u_{1},u_{2},v_{1},v_{2},\nu_{1},\nu_{2})\nonumber\\
&&=f(x)+\theta(y)-\langle\lambda,Ex+Fy+Gz-c\rangle+u_{1}^{\top}(x-u)-u_{2}^{\top}(x-l)\nonumber\\
&&~~~+v_{1}^{\top}(y-q)-v_{2}^{\top}(y-p)+\nu_{1}^{\top}(z-s)-\nu_{2}^{\top}(z-r),
\end{eqnarray}
where $(x,y,z)$  and $(\lambda=(\lambda^{e},\lambda^{ie}),u_{1},u_{2},v_{1},v_{2},\nu_{1},\nu_{2})$ are called the primal and dual variables.
%The relaxation augmented Lagrangian function of (\ref{y0}) is
%\begin{eqnarray}\label{y1}
%\mathcal{L}_{\beta}(x,y,z,\lambda)&=&f(x)+\theta(y)-\lambda^{\top}(Ex+Fy+Gz-c)+\frac{\beta}{2}\|Ex+Fy+Gz-c\|^{2}\nonumber\\
%&=&f(x)+\theta(y)+\frac{\beta}{2}\|Ex+Fy+Gz-c-\frac{\lambda}{\beta}\|^{2}-\frac{1}{2\beta}\|\lambda\|^{2},
%\end{eqnarray}
%where $(x,y,z)$  and $\lambda=((\lambda^{e})^T;(\lambda^{ie})^T)^T $ are called the primal and dual variables, and $\beta>0$ is a penalty parameter.
%Let $\begin {array}{ll} \lambda=\left(\ba{c} \lambda^{e}\\ \lambda^{ie} \ea\right),~\lambda_{k}=\left(\ba{c} \lambda^{e}_{k}\\ \lambda^{ie}_{k} \ea\right)\end {array}$ denote the Lagrangian multiplier. The full Lagrangian function of problem~(\ref{y0}) is defined as
%\begin{eqnarray}\label{1.3}
%&&{L}(x,y,z,\lambda,u_{1},u_{2},v_{1},v_{2},\nu_{1},\nu_{2})\nonumber\\
%&&=f(x)+\theta(y)-\langle\lambda,Ex+Fy+Gz-c\rangle+u_{1}^{\top}(x-u)-u_{2}^{\top}(x-l)\nonumber\\
%&&~~~+v_{1}^{\top}(y-q)-v_{2}^{\top}(y-p)+\nu_{1}^{\top}(z-s)-\nu_{2}^{\top}(z-r).
%\end{eqnarray}
Furthermore, for convenience of expression and analysis, we define operations between $\{-\infty,+\infty\}$ and the real number set $\mathfrak {R}$ as follows:
%In order to deal with the normal and generalized boundaries of box constraints uniformly under the optimality conditions and related analysis, the operations of~$\{\pm\infty\}$ and~$\mathfrak {R}^{n}$ are defined as follow:
$$\ba{cc}
\pm\infty+a=\pm\infty,~\forall~a\in\mathfrak {R};\pm\infty\times a=\pm\infty,~\forall~a>0;\\
\pm\infty\times a=\mp\infty,~\forall~a<0;\ \pm\infty\times a=0\Leftrightarrow~a=0.
\ea$$

For the current iterate~$(x_{k},y_{k},z_{k})$, where~$l\leq x_{k}\leq u, p\leq y_{k}\leq q, r\leq z_{k}\leq s,$ the traditional SQP method for the problem~(\ref{1.2}) is solving the following quadratic programming (QP)~subproblem:
\be\label{1.4}
\ba{ll} \min &\nabla f(x_{k})^{\top}(x-x_{k})+\nabla \theta(y_{k})^{\top}(y-y_{k})+\frac{1}{2}\|(x-x_{k},y-y_{k},z-z_{k})\|^{2}_{\mathcal{H}_{k}}\\
{\rm s.t.}&Ex+Fy+Gz=c,\
l\leq x\leq u,
p\leq y\leq q,
r\leq z\leq s,
\ea
\ee
where $\mathcal{H}_{k}$ is the symmetric approximation of the Hessian $$ \nabla^{2}_{(x,y,z)}{L}(\cdot)={\rm{diag}}(\nabla^{2}f(x),\nabla^{2}\theta(y),0)$$ for the full Lagrangian function (\ref{1.3}) with respect to ~$(x,y,z)$. Hence, the most preferable choice of the matrix~$\mathcal{H}_{k}$ is
     $\mathcal{H}_{k}={\rm{diag}}(H_{k}^{x},H_{k}^{y},0),$
     where~$H_{k} ^{x}$ and~$H_{k} ^{y}$ are the symmetric approximation matrices of $\nabla^{2}f(x_{k})$ and $\nabla^{2}\theta(y_{k})$, respectively.
As a result, the problem (\ref{1.4}) can be reduced as a three-block problem as follows:
\be\label{1.5}
\ba{ll} \min\ & \nabla f(x_{k})^{\top}(x-x_{k})+\frac{1}{2}\|x-x_{k}\|^{2}_{H_{k}^{x}}+\nabla\theta(y_{k})^{\top}(y-y_{k})+\frac{1}{2}\|y-y^{k}\|^{2}_{H_{k}^{y}}\\
{\rm s.t.~}\ &   Ex+Fy+Gz=c,\ l\leq x\leq u,
p\leq y\leq q,
r\leq z\leq s.
\ea
\ee
And the relaxed (ignoring the box constraints) augmented Lagrangian function of the QP subproblem (\ref{1.5}) is
\begin{eqnarray}\label{001.10}
&&\mathcal{L}_{\beta}^{\rm SQP}(x,y,z,\lambda)\nonumber\\
&=&\nabla f(x_{k})^{\top}(x-x_{k})+\frac{1}{2}\|x-x_{k}\|^{2}_{H_{k}^{x}}+\nabla \theta(y_{k})^{\top}(y-y_{k})+\frac{1}{2}\|y-y_{k}\|^{2}_{H_{k}^{y}}\nonumber\\
&&-\lambda^{\top}(Ex+Fy+Gz-c)+\frac{\beta}{2}\|Ex+Fy+Gz-c\|^{2}\nonumber\\
&=&\nabla f(x_{k})^{\top}(x-x_{k})+\frac{1}{2}\|x-x_{k}\|^{2}_{H_{k}^{x}}+\nabla \theta(y_{k})^{\top}(y-y_{k})+\frac{1}{2}\|y-y_{k}\|^{2}_{H_{k}^{y}}\nonumber\\
&&+\frac{\beta}{2}\|Ex+Fy+Gz-c-\frac{\lambda}{\beta}\|^{2}-\frac{1}{2\beta}\|\lambda\|^{2},
\end{eqnarray}
where $\beta>0$ is a penalty parameter.

%Based on~(\ref{001.10}) and the technology referred to in the reference \cite{Jian_Lao_Chao_Ma_2018}, we use ADMM and augmented Lagrangian thought to split (\ref{1.5}) into three subproblems as follows:

To reduce the calculation cost of the QP above, especially for large-scale problems,
on the basis of the decomposition idea of ADMM, we consider to split (\ref{1.5}) with respect to the variables $x,\ y$ and $z$ by the Jacobian method into three small-scale subproblems:
$$
\left\{\begin{array}{l}
\min\{\mathcal{L}_{\beta}^{\rm SQP}(x,y_{k},z_{k},\lambda_{k})|\ l\leq x\leq u\},\nonumber\\
\min\{\mathcal{L}_{\beta}^{\rm SQP}(x_{k},y,z_{k},\lambda_{k})|\ p\leq y\leq q\},\nonumber\\
\min\{\mathcal{L}_{\beta}^{\rm SQP}(x_{k},y_{k},z,\lambda_{k})|\ r\leq z\leq s\},\nonumber\\
\end{array}
    \right.
$$
where $\lambda_k=(\lambda_k^e,\lambda_k^{ie})\in \mathfrak {R}^{m_{1}}\times{R}^{m_{2}}$ is the multiplier vector corresponding to the equality and inequality constraints of the problem (2.1).
Obviously, the above three QP subproblems can be described as follows, respectively:
\be\label{1.6}
x-{\rm QP}\ba{ll} \min\ \nabla f(x_{k})^{\top}(x-x_{k})+\frac{1}{2}\|x-x_{k}\|^{2}_{H_{k}^{x}}+\frac{\beta}{2}\|Ex+Fy_{k}+Gz_{k}-c-\frac{\lambda_{k}}{\beta}\|^{2}\\\
{\rm s.t.~}\  l\leq x\leq u;
\ea
\ee
\be\label{1.7}
y-{\rm QP}\ba{ll} \min\ \nabla \theta(y_{k})^{\top}(y-y_{k})+\frac{1}{2}\|y-y_{k}\|^{2}_{H_{k}^{y}}+\frac{\beta}{2}\|Ex_{k}+Fy+Gz_{k}-c-\frac{\lambda_{k}}{\beta}\|^{2}\\\
{\rm s.t.~}\  p\leq y\leq q;
\ea
\ee
and
\be\label{1.8} z-{\rm QP}\ \  \ba{ll} \min\ \frac{\beta}{2}\|Ex_{k}+Fy_{k}+Gz-c-\frac{\lambda_{k}}{\beta}\|^{2}\\\
{\rm s.t.~}\  r\leq z\leq s.
\ea
\ee
Furthermore, after simple algebraic operations on the objective function of (\ref{1.8}), then the $z-{\rm QP}$ subproblem above can be reduced equivalently to
\be\label{1.11}\ba{ll} \min\{ \frac{\beta}{2}\|z-\hat{z}_{k+1}\|^{2}|\ r\leq z\leq s\}\\
\ea
\ee
with
\begin{equation}\label{x3}
\hat{z}_{k+1}=Cx_{k}+Dy_{k}-\frac{\lambda_{k}^{ie}}{\beta}.
\end{equation}
Next, in order to guarantee the solvability of the QP subproblems (\ref{1.6}) and (\ref{1.7}),
the matrices $H^{x}_k$ and $~H^{y}_k$ need to satsfiy~$(H^{x}_k+\beta E^{\top}E)\succ 0$ and $~(H^{y}_k+\beta F^{\top}F)\succ 0.$
As a result, the two QPs (\ref{1.6}) and (\ref{1.7}) are feasible and strictly convex, and each one has a unique optimal solution,
%
%the following assumptions should be made:
%\begin{assumption}\label{ass1} {\rm
%The matrices $H^{x}_k$ and $~H^{y}_k$ satsfiy~$(H^{x}_k+\beta E^{\top}E)\succ 0$ and $~(H^{y}_k+\beta F^{\top}F)\succ 0,$ i.e., ~$H^{x}_k+\beta(A^{\top}A+C^{\top}C)\succ 0,~H^{y}_k+\beta (B^{\top}B+D^{\top}D)\succ 0.$
%}
%\end{assumption}
%
%A sufficient condition for the assumption \ref{ass1} to hold is given below.
%
%{\bf Assumption~\ref{ass1}A.}\ \ (1)\ $H^k_x$ is a positive definite matrix or $E^{\top}E\succ0$ and ~$\beta>\frac{\|H^k_x\|}{\gamma_{\min}(E^{\top}E)}$;
%
% (2)\ $H^k_y$ is a positive definite matrix or $F^{\top}F\succ0$ and $\beta>\frac{\|H^k_y\|}{\gamma_{\min}(F^{\top}F)}$~, where the operator $\gamma_{\min}(\cdot)$ represents the minimum eigenvalue of the matrix.
%
%Obviously, both (\ref{1.6}) and (\ref{1.7}) have unique solution under the Assumption \ref{ass1}.
and let their optimal solutions be $\tilde{x}_{k+1}$ and $\tilde{y}_{k+1}$, respectively.  Again,
taking into account the special structure of the subproblem (\ref{1.11}), its optimal solution $\tilde{z}_{k+1}=\left((\tilde{z}_{k+1})_{i},~i=1,\ldots,m_{2}\right)$ can be generated by the following explicit closed expressions:
\be\label{3.9}(\tilde{z}_{k+1})_{i}=P_{[r_i,s_i]}((\hat{z}_{k+1})_i)=\left\{\begin{array}{lll}
(\hat{z}_{k+1})_{i},&\hbox{if}\ (\hat{z}_{k+1})_{i}\in[r_{i},s_{i}];\\
r_{i},&\hbox{if}\ (\hat{z}_{k+1})_{i}<r_{i};\\
s_{i},&\hbox{if}\ (\hat{z}_{k+1})_{i}>s_{i},
\end{array}
    \right.
    \ee
where $P_{[r_i,s_i]}(\cdot)$ is the projection operator.% from the point $z$ into the closed interval $[r,s]$.

Now,  according to the KKT optimality conditions for the subproblems (\ref {1.6})-(\ref{1.8}), there exist multiplier vectors $ \alpha^k_x,\gamma^k_x\in\mathfrak {R}^{n_{1}}$, $\alpha ^k_y,\gamma^k_y\in\mathfrak {R}^{n_{2}}$ and $ \alpha^k_z,\gamma^k_z\in\mathfrak {R}^{m_{2}}$ such that
\begin{subequations} \label{2.2}
\begin{numcases}{} \nabla f(x_k)+H^{x}_k(\tilde{x}_{k+1}-x_k)-E^{\top}[\lambda_k-\beta(E\tilde{x}_{k+1}+Fy_{k}+Gz_{k}-c)]-\nonumber\\
\alpha^{x}_{k}+\gamma^{x}_{k}=0,  \label{2.2a}\\
0\leq\alpha_{k}^{x}\perp(\tilde{x}_{k+1}-l)\geq 0,~0\leq\gamma_{k}^{x}\perp(u-\tilde{x}_{k+1})\geq 0; \label{2.2b}
\end{numcases}
\end{subequations}
\begin{subequations} \label{2.3}
\begin{numcases}{} \nabla \theta(y_k)+H^{y}_k(\tilde{y}_{k+1}-y_k)-F^{\top}[\lambda_k-\beta(Ex_{k}+F\tilde{y}_{k+1}+Gz_{k}-c)]-\nonumber\\\alpha^{y}_{k}+\gamma^{y}_{k}=0, \label{2.3a}\\
0\leq\alpha_{k}^{y}\perp(\tilde{y}_{k+1}-p)\geq 0,~0\leq\gamma_{k}^{y}\perp(q-\tilde{y}_{k+1})\geq 0;
 \label{2.3b}
\end{numcases}
\end{subequations}
and
\begin{subequations} \label{3.0}
\begin{numcases}{} \beta G^{\top}(Ex_{k}+Fy_{k}+G\tilde{z}_{k+1}-c-\frac{\lambda_{k}}{\beta})-\alpha^{z}_{k}+\gamma^{z}_{k}=0, \label{3.0a}\\
0\leq\alpha^{z}_{k}\perp(\tilde{z}_{k+1}-r)\geq0,~0\leq\gamma_{k}^{z}\perp(s-\tilde{z}_{k+1})\geq 0, \label{3.0b}
\end{numcases}
\end{subequations}
$$\Leftrightarrow\left\{ \ba{ll} \beta(\tilde{z}_{k+1}-\hat{z}_{k+1})-\alpha^{z}_{k}+\gamma^{z}_{k}=0,\\\
0\leq\alpha^{z}_{k}\perp(\tilde{z}_{k+1}-r)\geq0,~0\leq\gamma_{k}^{z}\perp(s-\tilde{z}_{k+1})\geq 0.
\ea
\right.
$$
On the basis of the KKT optimality conditions (\ref{2.2})--(\ref{3.0}), we analyze the descent property of the relaxation augmented Lagrangian function $\mathcal{L}_{\beta}(x,y,z,\lambda)$ of the problem (\ref{1.2}) along direction \be \label{2.4} d_{k}^{x}:=\tilde{x}_{k+1}-x_{k},~d_{k}^{y}:=\tilde{y}_{k+1}-y_{k},~d_{k}^{z}:=\tilde{z}_{k+1}-z_{k}
\ee at the current iteration point $(x_{k},y_{k},z_{k})$ , where $\mathcal{L}_{\beta}(x,y,z,\lambda)$ is defined as
\begin{eqnarray}\label{y1}
&\mathcal{L}_{\beta}&(x,y,z,\lambda)\nonumber\\&=&f(x)+\theta(y)-\lambda^{\top}(Ex+Fy+Gz-c)+\frac{\beta}{2}\|Ex+Fy+Gz-c\|^{2}\nonumber\\
&=&f(x)+\theta(y)+\frac{\beta}{2}\|Ex+Fy+Gz-c-\frac{\lambda}{\beta}\|^{2}-\frac{1}{2\beta}\|\lambda\|^{2}.
\end{eqnarray}
%So, the gradient of the relaxation augmented Lagrangian function $\mathcal{L}_{\beta}(x,y,z,\lambda)$ with respect to $(x,y,z)$
%\be \label{2.5}
%\nabla_{x}\mathcal{L}_{\beta}(x_{k},y_{k},z_{k},\lambda_{k})=\nabla f(x_{k})-E^{\top}[\lambda_{k}-\beta(Ex_{k}+Fy_{k}+Gz_{k}-c)],
%\ee
%\be \label{2.6}
%\nabla_{y}\mathcal{L}_{\beta}(x_{k},y_{k},z_{k},\lambda_{k})=\nabla \theta(y_{k})-F^{\top}[\lambda_{k}-\beta(Ex_{k}+Fy_{k}+Gz_{k}-c)],
%\ee
%\be \label{2.7}
%\nabla_{z}\mathcal{L}_{\beta}(x_{k},y_{k},z_{k},\lambda_{k})=-G^{\top}[\lambda_{k}-\beta(Ex_{k}+Fy_{k}+Gz_{k}-c)]
%=\beta(z_{k}-\hat{z}_{k+1}).
%\ee
From (\ref{y1}), (\ref{2.2a}),~(\ref{2.3a}), (\ref{3.0a}) and (\ref{2.4}), we can conclude that
\be \label{2.8}\ba{ll}
\nabla_{x}\mathcal{L}_{\beta}(x_{k},y_{k},z_{k},\lambda_{k})&
\overset{\eqref{y1}}{=}\nabla f(x_{k})-E^{\top}[\lambda_{k}-\beta(Ex_{k}+Fy_{k}+Gz_{k}-c)]\\&
\overset{\eqref{2.2a}}{=}
-H_{k}^{x}d_{k}^{x}-\beta E^{\top}Ed_{k}^{x}+\alpha_{k}^{x}-\gamma_{k}^{x},
\ea\ee
\be \label{2.9}\ba{ll}
\nabla_{y}\mathcal{L}_{\beta}(x_{k},y_{k},z_{k},\lambda_{k})&
\overset{\eqref{y1}}{=}
\nabla \theta(y_{k})-F^{\top}[\lambda_{k}-\beta(Ex_{k}+Fy_{k}+Gz_{k}-c)]\\&
\overset{\eqref{2.3a}}{=}
-H_{k}^{y}d_{k}^{y}-\beta F^{\top}Fd_{k}^{y}+\alpha_{k}^{y}-\gamma_{k}^{y},
\ea\ee
\begin{equation} \label{3.1}
\nabla_{z}\mathcal{L}_{\beta}(x_{k},y_{k},z_{k},\lambda_{k})
\overset{\eqref{y1}}{=}
-G^{\top}[\lambda_{k}-\beta(Ex_{k}+Fy_{k}+Gz_{k}-c)]
\overset{\eqref{3.0a}}{=}
-\beta d_{k}^{z}+\alpha_{k}^{z}-\gamma_{k}^{z}.
\end{equation}
Therefore, by (\ref{2.2b}), (\ref{2.4}) and (\ref{2.8}), it follows that
\begin{eqnarray} \label{3.2}
&\nabla_{x}&\mathcal{L}_{\beta}(x_{k},y_{k},z_{k},\lambda_{k})^{\top}d_{k}^{x}\nonumber\\&=&-\|d_{k}^{x}\|^{2}_{H_{k}^{x}}-\beta \|d_{k}^{x}\|^{2}_{E^{\top}E}-(u-x_{k})^{\top}\gamma^{x}_{k}+(l-x_{k})^{\top}\alpha_{k}^{x}\nonumber\\
&\leq&-\|d_{k}^{x}\|^{2}_{(H_{k}^{x}+\beta E^{\top}E)}.
\end{eqnarray}
Similarly, from (\ref{2.3b}), (\ref{2.4}) and (\ref{2.9}), one also gets
\begin{eqnarray} \label{3.3}
&\nabla_{y}&\mathcal{L}_{\beta}(x_{k},y_{k},z_{k},\lambda_{k})^{\top}d_{k}^{y}\nonumber\\
&=&-\|d_{k}^{y}\|^{2}_{H_{k}^{y}}-\beta \|d_{k}^{y}\|^{2}_{F^{\top}F}-(q-y_{k})^{\top}\gamma^{y}_{k}+(p-y_{k})^{\top}\alpha_{k}^{y}\nonumber\\
&\leq&-\|d_{k}^{y}\|^{2}_{(H_{k}^{y}+\beta F^{\top}F)}.
\end{eqnarray}
Obviously, it follows from (\ref{3.0b}), (\ref{2.4}) and (\ref{3.1}) that
\be\ba{ll} \label{3.4}
\nabla_{z}\mathcal{L}_{\beta}(x_{k},y_{k},z_{k},\lambda_{k})^{\top}d_{k}^{z}&=-\beta \|d_{k}^{z}\|^{2}-(s-z_{k})^{\top}\gamma_{k}^{z}+(r-z_{k})^{\top}\alpha_{k}^{z}\\&\leq-\beta\|d_{k}^{z}\|^{2}.\ea
\ee
Hence, it is easy to know from  (\ref{3.2})--(\ref{3.4}) that $\mathcal{L}_{\beta}(x,y,z,\lambda_k)$, with respect to
$(x,y,z)$, along the direction $(d_k^x,d_k^y,d_k^z)$ is descent at $(x_k, y_k,z_k)$.

In order to describe our analysis in a compact way,
denote
\begin{subequations} \label{3.5}
\begin{numcases}{}
u=(x,y,z),\ u_{k}=(x_{k},y_{k},z_{k}),
d_{k}^{u}=(d_{k}^{x},\ d_{k}^{y},d_{k}^{z}), \label{3.5a}\\
H_{k}^{u}={\rm diag}(
             H_{k}^{x}+\beta E^{\top}E,  H_{k}^{y}+\beta F^{\top}F, \beta I_{m_2}).\label{3.5b}
\end{numcases}
\end{subequations}
As a result, by (\ref{3.2})--(\ref{3.5}), it follows that
\be\label{3.6}
\nabla_{u}\mathcal{L}_{\beta}(u_{k},\lambda_{k})^{\top}d_{k}^{u}\leq
-\|d_{k}^{u}\|^{2}_{H_{k}^{u}}.
\ee
This inequality shows that $\mathcal{L}_{\beta}(\cdot,\lambda_k)$ along the direction $d_k^u$ has a nice descent property at $u_k$.
Hence, we consider to yield the next iteration point $x^{k+1}=(x_{k+1},y_{k+1},z_{k+1})$ by Armijo line search, according to the merit function $L_{\beta}(\cdot,\lambda^k)$ along direction $d_k^u$ at $u_{k}$.

%From assumption \ref{ass1} we can know that $H_{k}^{u}$ is a positive definite matrix, thus~$\mathcal{L}_{\beta}(u_{k},\lambda_{k})$ is descent along the direction $d_{k}^{u}$ at point $u_{k}$.

Subsequently, the other key problem to be addressed is how to yield $\lambda_{k+1}$ by updating $\lambda_{k}$. In this work, we consider the following correction:
%choose a suitable search direction for the relaxation augmented Lagrangian
%function (\ref{y1}) with respect to the dual variable $\lambda$. Considering that $\mathcal{L}_{\beta}(x,y,z,\lambda)$ is affine with respect to $\lambda$, we choose a reasonable deflection of the steepest descent direction $-\nabla_\lambda \mathcal{L}_{\beta}(x,y,z\lambda)$ as its improved direction for $\lambda$, i.e.,
\be\label{3.7}
\lambda_{k+1}=\lambda_k-\xi\nabla_{\lambda}\mathcal{L}_{\beta}(u_{k+1},\lambda_{k})=\lambda_{k}+\xi(Ex_{k+1}+Fy_{k+1}+Gz_{k+1}-c),
\ee
where the parameter~$\xi>0$, which is used to improve the numerical results.

To prepare for the subsequent analysis, we give the necessary conditions for KKT optimality of the problem (\ref{1.2}). A point $(\bar{x},\bar{y},\bar{z})$ is said to be a KKT point of the problem (\ref{1.2}) with a multiplier $\bar{\Lambda}:=(\bar{\lambda}:=(\bar{\lambda}^{e},\bar{\lambda}^{ie}), \bar{\alpha}^x, \bar{\gamma}^x,\bar{\alpha}^y, \bar{\gamma}^y,\bar{\alpha}^z, \bar{\gamma}^z)$,
%\in \mathbb{R}^m\times\mathbb{R}^{n_1}\times\mathbb{R}^{n_1}\times\mathbb{R}^{n_2}\times\mathbb{R}^{n_2}$,
if
\begin{equation}\label{2.3x}
\left\{
\begin{array}{l}
\nabla f(\bar{x})-E^{\top}\bar{\lambda}-\bar{\alpha}^x+\bar{\gamma}^x=0,\\
\nabla \theta(\bar{y})-F^{\top}\bar{\lambda}-\bar{\alpha}^y+\bar{\gamma}^y=0,\\ -G^{\top}\bar{\lambda}-\bar{\alpha}^z+\bar{\gamma}^z=0
,\\
0\leq \bar{\alpha}^x \perp (\bar{x}-l)\geq 0,\ 0\leq \bar{\gamma}^x \perp (u-\bar{x})\geq 0,\\
0\leq \bar{\alpha}^y \perp (\bar{y}-p)\geq 0,\ 0\leq \bar{\gamma}^y \perp (q-\bar{y})\geq 0,\\
0\leq \bar{\alpha}^z \perp (\bar{z}-r)\geq 0,\ 0\leq \bar{\gamma}^z \perp (s-\bar{z})\geq 0,\\
E\bar{x}+F\bar{y}+G\bar{z}=c.
\end{array}
\right.
\end{equation}
 Further, we call  $(\bar{x},\bar{y},\bar{z},\bar{\Lambda})$ satisying relationship \eqref{2.3x} a primal-dual solution to the problem (\ref{1.2}).
Based on the KKT optimization conditions above, the following lemma is at hand.
\begin{lemma}\label{lem1x}{\rm If $(\bar{x},\bar{y},\bar{z})$ is a KKT point of the problem (\ref{1.2}) with a multiplier $\bar{\Lambda}$,
%$$(\bar{\lambda}:=(\bar{\lambda}^e,\bar{\lambda}^{ie}), \bar{\alpha}^x, \bar{\gamma}^x,\bar{\alpha}^y, \bar{\gamma}^y,\bar{\alpha}^z, \bar{\gamma}^z),$$
 then $(\bar{x},\bar{y})$ is a KKT point of the problem (\ref{1.1}) with a multiplier $\bar{\Lambda}^e:=(\bar{\lambda}^e, \bar{\alpha}^x, \bar{\gamma}^x,\bar{\alpha}^y, \bar{\gamma}^y,\bar{\alpha}^z, \bar{\gamma}^z)$, namely,
\begin{equation}\label{2.3xb}
\left\{
\begin{array}{l}
\left(
\begin{array}{c}
\nabla f(\bar{x})\\
\nabla \theta(\bar{y})
\end{array}
\right)-
\left(
\begin{array}{c}
A^{\top}\\
B^{\top}
\end{array}
\right)\bar{\lambda}^e
+\left(
\begin{array}{c}
C^{\top}\\
D^{\top}
\end{array}
\right)
(\bar{\gamma}^z-\bar{\alpha}^z)+
\left(
\begin{array}{c}
\bar{\gamma}^x-\bar{\alpha}^x\\
\bar{\gamma}^y-\bar{\alpha}^y
\end{array}
\right)
=
\left(
\begin{array}{c}
0\\
0
\end{array}
\right),\\
0\leq \bar{\alpha}^x \perp (\bar{x}-l)\geq 0,\ 0\leq \bar{\gamma}^x \perp (u-\bar{x})\geq 0,\\
0\leq \bar{\alpha}^y \perp (\bar{y}-p)\geq 0,\ 0\leq \bar{\gamma}^y \perp (q-\bar{y})\geq 0,\\
0\leq \bar{\alpha}^z \perp (C\bar{x}+D\bar{y}-r)\geq 0,\ 0\leq \bar{\gamma}^z \perp (s-(C\bar{x}+D\bar{y}))\geq 0,\\
A\bar{x}+B\bar{y}=b.
\end{array}
\right.
\end{equation}
Further, $(\bar{x},\bar{y},\bar{\Lambda}^e)$ satisying relationship \eqref{2.3xb} is said to be a primal-dual solution to the problem (\ref{1.1}).
}
\end{lemma}

%We then have the following result, which plays an
%important role in our analysis:
\begin{lemma}\label{lem1}{\rm Suppose that $(H^{x}_k+\beta E^{\top}E)\succ 0$ and $~(H^{y}_k+\beta F^{\top}F)\succ 0.$
If the direction $d_{k}^{u}$ generated by (\ref{2.4})~and~(\ref{3.5})~equals zero, and $Ex_k+Fy_k+Gz_k-c=0$, then $(x_{k},y_k)$ is the KKT point of (\ref{1.1}) with the corresponding multiplier $(\lambda_{k}^e$,$\alpha_{k}^{x},$$\gamma_{k}^{x},$$\alpha_{k}^{y},$$\gamma_{k}^{y},$$\alpha_{k}^{z},$$\gamma_{k}^{z})$.}
\end{lemma}

{\bf Proof}~~From~(\ref{2.4}) and the stated conditions, one has
$$
\tilde{x}_{k+1}=x_{k},~~\tilde{y}_{k+1}=y_{k},~~\tilde{z}_{k+1}=z_{k},~~Ex_{k}+Fy_{k}+Gz_{k}=c.
$$
This shows that $u_{k}$ is a feasible solution to (\ref{1.2}).
Furthermore, taking into account the KKT conditions~(\ref{2.2})-(\ref{3.0}), one can gets
$$
\left\{\begin{array}{l}
\nabla f(x_{k})-E^{\top}\lambda_{k}-\alpha_{k}^{x}+\gamma_{k}^{x}=0,\nonumber\
\nabla \theta(y_{k})-F^{\top}\lambda_{k}-\alpha_{k}^{y}+\gamma_{k}^{y}=0,\nonumber\\
-G^{\top}\lambda_{k}-\alpha_{k}^{z}+\gamma_{k}^{z}=0,\nonumber\
0\leq\alpha_{k}^{x}\perp({x}_{k}-l)\geq 0,~0\leq\gamma_{k}^{x}\perp(u-{x}_{k})\geq 0,\nonumber\\
0\leq\alpha_{k}^{y}\perp({y}_{k}-p)\geq 0,~0\leq\gamma_{k}^{y}\perp(q-{y}_{k})\geq 0,\nonumber\\
0\leq\alpha_{k}^{z}\perp({z}_{k}-r)\geq 0,~0\leq\gamma_{k}^{z}\perp(s-{z}_{k})\geq 0.\nonumber
\end{array}
    \right.
$$
This, together with $Ex_k+Fy_k+Gz_k=c$ and (\ref{2.3x}), shows that $u_{k}$ is a KKT point of (\ref{1.2}) with the corresponding multiplier ~$(\lambda_{k},\alpha_{k}^{x},\gamma_{k}^{x},\alpha_{k}^{y},\gamma_{k}^{y},\alpha_{k}^{z},\gamma_{k}^{z}).$ Therefore, by Lemma \ref{lem1x}, the given conclusion is at hand.
\hfill$\Box$
%(ii)~From (\ref{3.7}), it is easy to get that
%\be\label{4.0}
%\nabla_{\lambda}\mathcal{L}_{\beta}(x_{k},y_{k},z_{k},\lambda_{k})^{\top}d_{k}^{\lambda}=-\frac{1}{\xi}\|d_{k}^{\lambda}\|^{2}.
%\ee
%Further, we denote\be\label{4.4}w=(u,\lambda),~w_{k}=(u_{k},\lambda_{k}),~H_{k}= \left(\begin{array}{ll}
%             H_{k}^{u}\ &  \\
%              \ & \frac{1}{\xi}I_{m_{1}+m_{2}}
%     \end{array}\right).
%\ee
%This, together with (\ref{3.6}) and (\ref{4.0}), implies that
%\begin{eqnarray}\label{4.3}
%\nabla\mathcal{L}_{\beta}(w_{k})^{\top}d_{k}\leq-\|d_{k}\|^{2}_{H_{k}}.
%\end{eqnarray}
%Combining the assumptions of $(H^{x}_k+\beta E^{\top}E)\succ 0$ and $~(H^{y}_k+\beta F^{\top}F)\succ 0$, so it is follows that $d_{k}$ is a descent direction of~$\mathcal{L}_{\beta}(w)$ at $w_{k}$. And the whole proof is completed.\hfill$\Box$

Now, on the basis of the analysis above, we give the steps of our splitting SQP (S-SQP) algorithm for solving two-block nonconvex optimization with general linear constraints (GLC) (GLC-S-SQP algorithm for short) as follows.

{\bf GLC-S-SQP algorithm A}

\indent{\bf{Step 0}}~(Initialization) Choose parameters~$\rho,\sigma\in(0,1),~\beta,\xi>0$ and initial iteration point $w_{0}:=(u_0,\lambda_{0})=(x_{0},y_{0},z_{0},\lambda_{0}^{e},\lambda_{0}^{ie})$ satisfying:~$l\leq x_{0}\leq u,~p\leq y_{0}\leq q,~r\leq z_{0}\leq s,$ two symmetric matrices $H_{0}^{x}\in \mathfrak {R}^{n_{1}\times n_{1}}$ and $H_{0}^{y}\in \mathfrak {R}^{n_{2}\times n_{2}}$ such that~$H_{0}^{x}+\beta E^{\top}E\succ 0,~H_{0}^{y}+\beta F^{\top}F\succ 0.$ Set~$k=0.$

{\bf{Step 1}}~(Solving~QPs)
Solving the two~${\rm QP}$ subproblems~(\ref{1.6}) and~(\ref{1.7}) in parallel to generate the (unique) optimal solutions $\tilde{x}_{k+1}$ and $\tilde{y}_{k+1},$ and yield $\tilde{z}_{k+1}$ by (\ref{x3})--(\ref{3.9}).

{\bf{Step 2}}~(Computing search direction)
Generate search direction $d_{k}^{u}$ by (\ref{3.5}).
%i.e.,
%\[d_{k}^{x}=\tilde{x}_{k+1}-{x}_{k},~d_{k}^{y}=\tilde{y}_{k+1}-{y}_{k},~d_{k}^{z}=\tilde{z}_{k+1}-{z}_{k},\]
%\[d_{k}^{\lambda}=\xi(E{x}_{k}+F{y}_{k}+G{z}_{k}-c),~d_{k}^{u}=(d_{k}^{x},d_{k}^{y},d_{k}^{z}),~\]
If~$d_{k}^u=0$ and $Ex_k+Fy_k+Gz_k-c=0$, then~$(x_{k},y_k)$ is a KKT point of (\ref{1.1}), stop. Otherwise, go to Step 3.

{\bf{Step 3}}~(Computing the step size) Compute the step size $t_k$ by Armijo line search, that is, the maximum $t$ of the sequence $\{1,\sigma,\sigma^{2},\ldots\}$ satisfying
\be\label{4.2}
\mathcal{L}_{\beta}(u_{k}+td^{u}_{k},\lambda_{k})\leq \mathcal{L}_{\beta}(u_{k},\lambda_{k})-t\rho\|d^{u}_{k}\|^{2}_{H^{u}_{k}}.
\ee

{\bf{Step 4}}~(Updating) Set $w_{k+1}=(u_{k+1},\lambda_{k+1})$ with $u_{k+1}=u_{k}+t_{k}d_{k}^u$ and
\be\label{2.30x}
\ba{ll}
\lambda_{k+1}=\left(\begin{array}{l}\lambda_{k+1}^e\\ \lambda_{k+1}^{ie}\end{array}\right)
&
=\lambda_k+\xi(Ex_{k+1}+Fy_{k+1}+Gz_{k+1}-c)\\
& =\left(\begin{array}{c}\lambda_{k}^e+\xi(Ax_{k+1}+By_{k+1}-b)\\ \lambda_{k}^{ie}+\xi(Cx_{k+1}+Dy_{k+1}-z_{k+1})\end{array}\right).
\ea
\ee
%Generate new iteration point, i.e.,
%$x_{k+1}=x_{k}+t_{k}d_{k}^{x},~y_{k+1}=y_{k}+t_{k}d_{k}^{y},~z_{k+1}=z_{k}+t_{k}d_{k}^{z},~\lambda_{k+1}=\lambda_{k}+t_{k}d_{k}^{\lambda},~u_{k+1}=(x_{k+1},y_{k+1},z_{k+1}),~w_{k+1}=(u_{k+1},\lambda_{k+1});$
Then generate new symmetric matrices $H_{k+1}^{x},~H_{k+1}^{y}$ such that~$H_{k+1}^{x}+\beta E^{\top}E\succ 0, H_{k+1}^{y}+\beta F^{\top}F\succ 0.$
Let~$k:=k+1,$ and go back to Step 1.

\begin{remark} {\rm In Step 3, if the direction $d_k^u=0$, the step length $t_k=1$, in this case, the primal variable $(x,y,z)$ is not updated,  and the dual primal variable $\lambda$ can be updated by \eqref{2.30x} in Step 4. Otherwise, by (\ref{3.6}), it is not difficult to know that (\ref{4.2}) is satisfied for sufficiently small $t>0$. So, the GLC-S-SQP algorithm A is well-defined.}
\end{remark}

%
%The following lemma illustrates that the line search in Step 3 is realizable.
%\begin{lemma}{\rm
%In the GLC-ADMM-SQP Algorithm~{\rm I}, (\ref{4.2}) hold so long as the positive number $t$ is sufficiently small. Thus, Step 3 can be terminated after a limited iteration, and the algorithm is well-defined.
%}
%\end{lemma}
%
%{\bf Proof}
%From Theorem~\ref{the333}, we know
%$$\mathcal{L}_{\beta}(w_{k}+td_{k})=\mathcal{L}_{\beta}(w_{k})+t\nabla\mathcal{L}_{\beta}(w_{k})^{\top}d_{k}+o(t\|d_{k}\|).$$
%For fixed $k$, we obtain that the following relation hold if $t$ is sufficiently small by combining above formula with (\ref{4.2})
%$$\mathcal{L}_{\beta}(w_{k}+td_{k})-\mathcal{L}_{\beta}(w_{k})+t\rho\nabla\mathcal{L}_{\beta}(w_{k})^{\top}d_{k}\leq-t(1-\rho)\|d_{k}\|_{H_{k}}^{2}+o(t).$$
%This completes the proof.\hfill$\Box$
%\setcounter{equation}{0}
\section{Convergence analysis}

If the GLC-S-SQP algorithm A stops at $w_k$, from Step 2 and Lemma \ref{lem1}, we know that the current iteration point $(x_k,y_k)$ is a KKT point of problem (\ref{1.1}). In this section, we assume that the algorithm yields an infinite iteration sequence $\{w_k\}$ of points, and analyze the global convergence of the GLC-S-SQP algorithm A under the following  basic assumption:

%From Lemma~\ref{lem1}, it is follows that the current iteration point~$u_{k}$ is the KKT~point of the problem (\ref{y0}) if GLC-ADMM-SQP algorithm {\rm I} terminates within a finite number of iterations. Therefore, we only need to prove the global convergence on the premise that the GLC-ADMM-SQP algorithm {\rm I} generates an infinite point sequence  $\{w_{k}\}$ in the following analysis.

\begin{assumption}\label{ass2}{\rm{ For any bounded subsequence $\{w_k\}_{\mathcal{K}}$ of $\{w_k\}$,
the associated sequences~$\{H_{k}^{x}\}_{\mathcal{K}}$ and~$\{H_{k}^{y}\}_{\mathcal{K}}$ of matrices in the GLC-S-SQP algorithm A are bounded, and there exist constants $\eta^{x}>0$ and $\eta^{y}>0$ such that}}
$$
H_{k}^{x}+\beta E^{\top}E\succ\eta^{x}I_{n_{1}},~H_{k}^{y}+\beta F^{\top}F\succ\eta^{y}I_{n_{2}},\ \forall\  k\in {\mathcal{K}}.
$$
\end{assumption}
It follows from~(\ref{3.5}) and Assumption~\ref{ass2} that
\be\label{4.8}
\|d_{k}^{u}\|^{2}_{H_{k}^{u}}\geq\tilde{\eta}\|d_{k}^{u}\|^{2},\ \forall\  k\in {\mathcal{K}}, \hbox{\ with}
~\tilde{\eta}=\min\{\eta^{x},\eta^{y},\beta\}>0.
\ee

\begin{lemma}\label{lem2}{\rm
Suppose that~Assumption \ref{ass2} holds. If a subsequence $\{w_{k}\}_{\mathcal{K}}$ of $\{w_{k}\}$ generated by the GLC-S-SQP algorithm A is bounded, then the corresponding subsequences $\{d_{k}^u\}_{\mathcal{K}},$
$\{w_{k+1}\}_{\mathcal{K}},\{\tilde{u}_{k+1}:=(\tilde{x}_{k+1},\tilde{y}_{k+1},\tilde{z}_{k+1})\}_{\mathcal{K}}$ and $\{(\alpha_{k}^{x},\gamma_{k}^{x},\alpha_{k}^{y},\gamma_{k}^{y},\alpha_{k}^{z},\gamma_{k}^{z})\}_{\mathcal{K}}$ are all bounded.
}
\end{lemma}

{\bf Proof}\ We first prove the boundedness of the sequence $\{d_{k}^{y}\}_{\mathcal{K}}$.
In view of $\tilde{y}_{k+1}$ and ${y}_{k}$ are optimal and feasible solutions to the subproblem (\ref{1.7}), respectively, we can obtain
\begin{eqnarray}\nonumber
&\nabla \theta(y_{k})^{\top}d_{k}^{y}+\frac{1}{2}\|d_{k}^{y}\|^{2}_{H_{k}^{y}}+\frac{\beta}{2}\|Ex_{k}+F\tilde{y}_{k+1}+Gz_{k}-c-\frac{\lambda_{k}}{\beta}\|^{2}\\&\leq\frac{\beta}{2}\|Ex_{k}+Fy_{k}+Gz_{k}-c-\frac{\lambda_{k}}{\beta}\|^{2}.\nonumber
\end{eqnarray}
Denote that~$a_{k}=\frac{\beta}{2}\|Ex_{k}+Fy_{k}+Gz_{k}-c-\frac{\lambda_{k}}{\beta}\|^{2}.$
Since the boundedness of $\{w_{k}\}_{\mathcal{K}}$, there exists a constant $M>0$ such that
\[\|\nabla \theta(y_{k})\|\leq M,~a_{k}\leq M,~ \|F^{\top}(Ex_{k}+Fy_{k}+Gz_{k}-c-\frac{\lambda_{k}}{\beta})\|\leq M,\ \forall\  k\in\mathcal{K}.\]
Therefore, for each~$k\in\mathcal{K}$, we have
\begin{eqnarray}\nonumber
M&\geq&\nabla \theta(y_{k})^{\top}d_{k}^{y}+\frac{1}{2}\|d_{k}^{y}\|^{2}_{H_{k}^{y}}+\frac{\beta}{2}\|Fd_{k}^{y}+Ex_{k}+Fy_{k}+Gz_{k}-c-\frac{\lambda_{k}}{\beta}\|^{2}\nonumber\\
&\geq&-\|\nabla \theta(y_{k})\|\cdot\|d_{k}^{y}\|+\frac{1}{2}\|d_{k}^{y}\|^{2}_{\left(H_{k}^{y}+\beta F^{\top}F\right)}-\beta\|d_{k}^{y}\|\nonumber\\&&\cdot\|F^{\top}(Ex_{k}+Fy_{k}+Gz_{k}-c-\frac{\lambda_{k}}{\beta})\|\nonumber\\
&&+\frac{\beta}{2}\|Ex_{k}+Fy_{k}+Gz_{k}-c-\frac{\lambda_{k}}{\beta}\|^{2}\nonumber\\
&\geq&-M\|d_{k}^{y}\|+\frac{1}{2}\|d_{k}^{y}\|^{2}_{\left(H_{k}^{y}+\beta F^{\top}F\right)}-\beta M\|d_{k}^{y}\|\nonumber\\
&=&\frac{1}{2}\|d_{k}^{y}\|^{2}_{\left(H_{k}^{y}+\beta F^{\top}F\right)}-(1+\beta)M\|d_{k}^{y}\|.\nonumber
\end{eqnarray}
%Further,~$\|d_{k}^{y}\|^{2}_{\left(H_{k}^{y}+\beta F^{\top}F\right)}-2(1+\beta)M\|d_{k}^{y}\|\leq2M.$
This, together with Assumption \ref{ass2},~ implies
\[\eta^{y}\|d_{k}^{y}\|^{2}-2(1+\beta)M\|d_{k}^{y}\|\leq2M, \forall\ k\in \mathcal{K}. \]
Therefore,~the boundedness of sequence $\{d_{k}^{y}\}_{\mathcal{K}}$ is at hand.
In a fashion similar to the analysis above,~the boundedness of $\{d_{k}^{x}\}_{\mathcal{K}}$ can be also proved.

On the other hand, it is easy to get that $\{\tilde{z}_{k+1}\}_{\mathcal{K}}$ is bounded from (\ref{x3}), (\ref{3.9}) and the boundedness of $\{w_{k}\}_{\mathcal{K}}$. Therefore,~$\{d_{k}^{z}\}_{\mathcal{K}}$ is also bounded. Further,~$\{d_{k}^u\}_{\mathcal{K}}$ is bounded.
Again,~$\{u_{k+1}=u_{k}+t_{k}d_{k}^{u}\}_{\mathcal{K}}$ and $\{\tilde{u}_{k+1}=u_{k}+d_{k}^{u}\}_{\mathcal{K}}$ are also bounded.
Furthermore, from (\ref{2.30x}) and the boundedness of $\{u_{k+1}\}_{\mathcal{K}}$ and $\{\lambda_{k}\}_{\mathcal{K}}$, the boundedness of $\{w_{k+1}\}_{\mathcal{K}}$ is at hand.

Finally, by~KKT condition~(\ref{2.3}) and the boundedness of $\{(w_{k},\tilde{u}_{k+1},d_k^u,H_{k}^{y})\}_{\mathcal{K}}$,
it follows that $\{(\alpha_{k}^{y}-\gamma_{k}^{y})\}_{\mathcal{K}}$ is bounded and ~$(\alpha_{k}^{y})^{\top}\gamma_{k}^{y}=0$.
Therefore, $\|\alpha_{k}^{y}\|^{2}=(\alpha_{k}^{y})^{\top}(\alpha_{k}^{y}-\gamma_{k}^{y})\leq \|\alpha_{k}^{y}\|\cdot\|\alpha_{k}^{y}-\gamma_{k}^{y}\|,$ which implies that $\{\alpha_{k}^{y}\}_{\mathcal{K}}$ is bounded,
 and so is $\{\gamma_{k}^{y}\}_{\mathcal{K}}$.
In a fashion similar to the analysis above,~sequences $\{(\alpha_{k}^{x},\gamma_{k}^{x})\}_{\mathcal{K}}$ and $\{(\alpha_{k}^{z},\gamma_{k}^{z})\}_{\mathcal{K}}$ are also bounded.  The whole proof is completed.\hfill$\Box$

%$\{(\alpha_{k}^{x}-\gamma_{k}^{x})\}_{\mathcal{K}}, \{(\alpha_{k}^{y}-\gamma_{k}^{y})\}_{\mathcal{K}}, \{(\alpha_{k}^{z}-\gamma_{k}^{z})\}_{\mathcal{K}}$ is bounded and ~$\alpha_{k}^{x}\perp\gamma_{k}^{x},~\alpha_{k}^{y}\perp\gamma_{k}^{y},~\alpha_{k}^{z}\perp\gamma_{k}^{z}$.

%Similarly, $\{\alpha_{k}^{x},\gamma_{k}^{x}\}_{\mathcal{K}}, \{\alpha_{k}^{z},\gamma_{k}^{z}\}_{\mathcal{K}}$ is also bounded. Therefore, $\{\gamma_{k}^{x}\}_{\mathcal{K}},~\{\gamma_{k}^{x}\}_{\mathcal{K}}$ and $~\{\gamma_{k}^{z}\}_{\mathcal{K}}$ is bounded. This completes the proof.\hfill$\Box$

The following analysis shows that the sequence $\{\mathcal{L}_{\beta}(w_{k})\}$ generated by the GLC-S-SQP algorithm A has nice monotonicity.
Taking into account the definition of $\mathcal{L}_{\beta}(\cdot)$, i.e, (\ref{y1}), we obtain, for $\forall\ x,y,z,\lambda$,
$$\mathcal{L}_{\beta}(x,y,z,\lambda+\xi(Ex+Fy+Gz-c))=\mathcal{L}_{\beta}(x,y,z,\lambda)-\xi\|Ex+Fy+Gz-c\|^{2}.$$
This, along with (\ref{4.2}) and \eqref{2.30x}, shows that
\be\label{3.55}
\mathcal{L}_{\beta}(w_{k+1})-\mathcal{L}_{\beta}(w_{k})\leq-\xi\|Ex_{k+1}+Fy_{k+1}+Gz_{k+1}-c\|^{2}-t_{k}\rho\|d^{u}_{k}\|^{2}_{H^{u}_{k}}~,\forall\ k\geq0.
\ee

Subsequently, we always assume that $w_{*}:=(x_{*},y_{*},z_{*},\lambda_{*})$ is a given accumulation point of the sequence $\{w_{k}\}$, then there exists an infinite subsequence $\mathcal{K}$ such that $w_{k}\rightarrow w_{*}, k\in\mathcal{K}.$ Therefore, from Lemma \ref{lem2}, we can assume, without loss of generality, that the following relations hold for $k\rightarrow\infty$ and $k\in\mathcal{K}$: %(otherwise, there exists a subsequence of $\mathcal{K}$)
\be\label{114.1}\left\{ \ba{ll}
w_{k}\rightarrow w_{*},~{d}_{k}^{u}:=(d_{k}^{x},d_{k}^{y},d_{k}^{z})\rightarrow {d}_{*}^{u}:=(d_{*}^{x},d_{*}^{y},d_{*}^{z}),\\
(\alpha_{k}^{x},\gamma_{k}^{x},\alpha_{k}^{y},\gamma_{k}^{y},\alpha_{k}^{z},\gamma_{k}^{z})\rightarrow(\alpha_{*}^{x},\gamma_{*}^{x}, \alpha_{*}^{y},\gamma_{*}^{y},\alpha_{*}^{z},\gamma_{*}^{z}).
\ea\right.
\ee

\begin{lemma}\label{lemma3.4}\rm{
Suppose that Assumption~\ref{ass2} holds, and the infinite index $\mathcal{K}$ such that relationship (\ref{114.1}) holds. Then the limit ${d}_{*}^{u}$ defined by (\ref{114.1}) equals zero, $w_{k+1}\rightarrow w_{*}, k\in\mathcal{K}$ and $Ex_{*}+Fy_{*}+Gz_{*}=c$.
}\end{lemma}

{\bf Proof}
Since $f$ and $\theta$ are continuously differentiable functions, it follows that $\mathcal{L}_{\beta}(w_{k})\rightarrow \mathcal{L}_{\beta}(w_{*}), k\in\mathcal{K}.$ This, together with (\ref{3.55}), implies that sequence $\{\mathcal{L}_{\beta}(w_{k})\}$
is monotonically decreasing and contains a convergent subsequence. Hence, the whole sequence $\{\mathcal{L}_{\beta}(w_{k})\}$ is convergent, furthermore, we have
$$\lim\limits_{k\in\mathcal{K}}\left(\mathcal{L}_{\beta}(w_{k+1})-\mathcal{L}_{\beta}(w_{k})\right)=0.$$
Passing to the limit in the inequality (\ref{3.55}) for $k\in\mathcal{K}$, it follows that
\be\ba{ll}\label{3.56}
0&=\lim\limits_{k\in\mathcal{K}}\left(\mathcal{L}_{\beta}(w_{k+1})-\mathcal{L}_{\beta}(w_{k})\right)\\&\leq\lim\limits_{k\in\mathcal{K}}\left(-\xi\|Ex_{k+1}+Fy_{k+1}+Gz_{k+1}-c\|^{2}-t_{k}\rho\|d^{u}_{k}\|^{2}_{H^{u}_{k}}\right).
\ea\ee
This, together
with~(\ref{4.8}), further shows that \be\label{3.561}\lim\limits_{k\in\mathcal{K}}t_{k}d^{u}_{k}=0, \lim\limits_{k\in\mathcal{K}}(Ex_{k+1}+Fy_{k+1}+Gz_{k+1}-c)=0.\ee

Next, we prove that $d_{*}^{u}=0.$ Assume that $d_{*}^{u}\neq0$ by contradiction. Then there exists $\varepsilon>0$ and $k_{0}\in\mathcal{K}$ such that $\|d^{u}_{k}\|>\varepsilon,~\forall\  k\in\mathcal{K}_{0}:=\{k\mid k\in\mathcal{K},~k>k_{0}\}.$ For $k\in\mathcal{K}_{0}$,
it follows, from Taylor expansion, (\ref{3.6}), (\ref{4.8}) and the boundedness of $\{d_{k}^{u}\}_{\mathcal{K}_{0}}$, that (for sufficiently small $t>0$ independent $k$)
\begin{eqnarray}\nonumber
\mathcal{L}_{\beta}(u_{k}+td_{k},\lambda_{k})&=&\mathcal{L}_{\beta}(u_{k},\lambda_{k})+t\nabla_{u} \mathcal{L}_{\beta}(u_{k},\lambda_{k})^{\top}d_{k}^{u}+o(t\|d_{k}^{u}\|)\nonumber\\
&\leq&\mathcal{L}_{\beta}(u_{k},\lambda_{k})-t\|d_{k}^{u}\|^{2}_{H_{k}^{u}}+o(t)\nonumber\\
&=&\mathcal{L}_{\beta}(u_{k},\lambda_{k})-t\rho\|d_{k}^{u}\|^{2}_{H_{k}^{u}}-t(1-\rho)\|d_{k}^{u}\|^{2}_{H_{k}^{u}}+o(t)\nonumber\\
&\leq&\mathcal{L}_{\beta}(u_{k},\lambda_{k})-t\rho\|d_{k}^{u}\|^{2}_{H_{k}^{u}}-t\tilde{\eta}(1-\rho)\|d_{k}^{u}\|^{2}+o(t)\nonumber\\
&\leq&\mathcal{L}_{\beta}(u_{k},\lambda_{k})-t\rho\|d_{k}^{u}\|^{2}_{H_{k}^{u}}-t\tilde{\eta}(1-\rho)\varepsilon^{2}+o(t)\nonumber\\
&\leq&\mathcal{L}_{\beta}(u_{k},\lambda_{k})-t\rho\|d_{k}^{u}\|^{2}_{H_{k}^{u}}.\nonumber
\end{eqnarray}
This, together with the Armijo line search rule (\ref{4.2}), implies that $t_{*}:=\inf\{t_{k}: k\in\mathcal{K}_{0}\}>0$. Hence, $\lim\limits_{k\in\mathcal{K}_{0}}\|t_{k}d_{k}^{u}\|\geq t_{*}\varepsilon>0,$ which contradicts the first relation of (\ref{3.561}).
So $d^{u}_{*}=0$ is at hand.
On the other hand, if follows from second relations of \eqref{2.30x} and (\ref{3.561}) that $\lim\limits_{k\in\mathcal{K}}\lambda_{k+1}=\lim\limits_{k\in\mathcal{K}}(\lambda_k+ \xi(Ex_{k+1}+Fy_{k+1}+Gz_{k+1}-c))=\lim\limits_{k\in\mathcal{K}}\lambda_{k}=\lambda_{*}$.
This, together with $(w_k,d_k^u){\overset{k\in\mathcal{K}}\longrightarrow}(w_*,0)$, shows that $w_{k+1}{\overset{k\in\mathcal{K}}\longrightarrow}w_{*}$.
Furthermore, this, along with the second relation of (\ref{3.561}) gives that $Ex_{*}+Fy_{*}+Gz_{*}=c$. And the whole proof is completed.\hfill$\Box$

Now, on the basis of Lemma \ref{lemma3.4}, we can obtain the global convergence of the GLC-S-SQP algorithm A as follows.
\begin{theorem}\rm{
Suppose that Assumption \ref{ass2} holds. Then for each accumulation point $w_{*}:=(x_*,y_*,z_*,\lambda_*=(\lambda_*^e, \lambda_*^{ie}))$ of the sequence $\{w_{k}\}$, $(x_*,y_*)$ is a KKT point of the problem (\ref{1.1}), and there exists an infinite
subsequence $\{(x_{k}, y_{k}, \lambda_{k}^e,\alpha_{k}^{x},\gamma_{k}^{x},\alpha_{k}^{y},\gamma_{k}^{y},\alpha_{k}^{z},\gamma_{k}^{z})\}_{\mathcal{K}}$ converges the primal-dual solution \\$(x_{*}, y_{*}, \lambda_{*}^e, \alpha_{*}^{x},\gamma_{*}^{x},\alpha_{*}^{y},\gamma_{*}^{y},\alpha_{*}^{z},\gamma_{*}^{z})$ to the problem (\ref{1.1}), i.e., the GLC-S-SQP algorithm A is globally convergent.
}\end{theorem}

{\bf Proof}  First, it follows from Lemma \ref{lemma3.4} that ${d}_{*}^u=0$ and $Ex_{*}+Fy_{*}+Gz_{*}=c$. Then, we obtain
\be\label{1.10}\left\{\ba{ll}\lim\limits_{k\in\mathcal{K}}\tilde{x}_{k+1}=\lim\limits_{k\in\mathcal{K}}(x_{k}+d_{k}^{x})=x_{*},
\lim\limits_{k\in\mathcal{K}}\tilde{y}_{k+1}=\lim\limits_{k\in\mathcal{K}}(y_{k}+d_{k}^{y})=y_{*},\\
\lim\limits_{k\in\mathcal{K}}\tilde{z}_{k+1}=\lim\limits_{k\in\mathcal{K}}(z_{k}+d_{k}^{z})=z_{*}.\ea\right.
\ee
These, together with $(x_\ast,y_\ast,z_\ast)\in [l,u]\times [p,q]\times [r,s]$, imply that $(x_*,y_*,z_*)$ is a feasible solution for  (\ref{1.2}). Next, passing to the limit in the KKT conditions~(\ref{2.2})-(\ref{3.0}) for $k\in\mathcal{K}$, respectively, we have
$$
\left\{\begin{array}{l}\nonumber
\nabla f(x_{*})-E^{\top}\lambda_*-\alpha_{*}^{x}+\gamma_{*}^{x}=0,\\
\nabla \theta(y_*)-F^{\top}\lambda_*-\alpha_{*}^{y}+\gamma_{*}^{y}=0,\\
-G^{\top}\lambda_{*}-\alpha_{*}^{z}+\gamma_{*}^{z}=0,\\
0\leq\alpha_{*}^{x}\perp({x}_{*}-l)\geq 0,~0\leq\gamma_{*}^{x}\perp(u-{x}_{*})\geq 0,\\
0\leq\alpha_{*}^{y}\perp({y}_{*}-p)\geq 0,~0\leq\gamma_{*}^{y}\perp(q-{y}_{*})\geq 0,\\
0\leq\alpha_{*}^{z}\perp({z}_{*}-r)\geq 0,~0\leq\gamma_{*}^{z}\perp(s-{z}_{*})\geq 0,\\
Ex_{*}+Fy_{*}+Gz_{*}-c=0.
\end{array} \right.
$$
These show that $(x_{*},y_{*},z_{*})$ with the multiplier $(\lambda_{*},\alpha_{*}^{x},\gamma_{*}^{x},\alpha_{*}^{y},\gamma_{*}^{y},\alpha_{*}^{z},\gamma_{k}^{z})$ satisfies (\ref{2.3x}). Moreover, by Lemma \ref{lem1x}, one knows that $(x_*,y_*)$ is a KKT point of the problem (\ref{1.1}), and $(x_{*}, y_{*}, \lambda_{*}^e, \alpha_{*}^{x},\gamma_{*}^{x},\alpha_{*}^{y},\gamma_{*}^{y},
\alpha_{*}^{z},\gamma_{*}^{z})$ is a primal-dual solution to the problem (\ref{1.1}).  The proof is completed. \hfill$\Box$

\section{An extension of the GLC-S-SQP algorithm A}
\label{sec:2}

In this section, the GLC-S-SQP algorithm A is further extended to solve the following optimization problem:
\be\label{111.1}
\min \{f(x)+\theta(y)\ | \ Ax+By=b,\ r\leq Cx+Dy\leq s,\
x\in \mathcal{X}, y\in \mathcal{Y}\},
\ee
where $\mathcal{X}\subseteq\mathfrak{R}^{n_{1}}$~and~$\mathcal{Y}\subseteq\mathfrak{R}^{n_{2}}$ are general nonempty closed convex sets.
%
%Similar to the analysis in section 3.1, the above problem is equivalent to:
%\be\label{111.2}
%\ba{ll} \min &f(x)+\theta(y)\\\
%\ {\rm s.t.}&Ex+Fy+Gz=c,\\
%&x\in \mathcal{X},~y\in \mathcal{Y},~z\geq 0.
%\ea
%\ee

Similar to the analysis in Section 2, the problem (\ref{111.1}) can be also reformulated as follows:
\be\label{111.1x}
\min \{f(x)+\theta(y)\ | \ Ex+Fy+Gz=c,\
x\in \mathcal{X}, y\in \mathcal{Y},\ r\leq z\leq s\}.
\ee

Then, from \cite{Rockafellar_Wets_1998}, the necessary optimality conditions of the problem (\ref{111.1x}) is given below.
% Then In a fashion similar to the analysis above, we give the optimality conditions of the problem (\ref{111.1x}).
A point $(\bar{x},\bar{y},\bar{z})$ is said to be a stationary point of the problem (\ref{111.1x}) with a multiplier $(\bar{\lambda}:=(\bar{\lambda}^{e},\bar{\lambda}^{ie}), \bar{\alpha}^z, \bar{\gamma}^z)$, if
\be\label{4.2x}
\left\{\begin{array}{ll}
0\in\nabla f(\bar{x})- E^{\top}\bar{\lambda}+ N_{\mathcal{X}}(\bar{x}),\\
0\in \nabla \theta(\bar{y})- F^{\top}\bar{\lambda}+ N_{\mathcal{Y}}(\bar{{y}}),\\
-G^{\top}\bar{\lambda}-\bar{\alpha}^{z}+\bar{\gamma}^{z}=0,\\
0\leq\bar{\alpha}^{z}\perp(\bar{{z}}-r)\geq 0,~0\leq\bar{\gamma}^{z}\perp(s-\bar{{z}})\geq 0,\\
E\bar{x}+F\bar{y}+G\bar{z}-c=0,
\end{array} \right.
\ee
where~$N_{\mathcal{X}}(\bar{x})$ and $N_{\mathcal{Y}}(\bar{y})$
are the normal cones of closed convex sets $\mathcal{X}$ and $\mathcal{Y}$ at the points $\bar{x}$ and $\bar{y}$, respectively.

Based on the optimality conditions above, the following lemma is at hand.
\begin{lemma}\label{lem1xxx}{\rm If $(\bar{x},\bar{y},\bar{z})$ is a stationary point of the problem (\ref{111.1x}) with a multiplier $(\bar{\lambda}:=(\bar{\lambda}^e,\bar{\lambda}^{ie}), \bar{\alpha}^z, \bar{\gamma}^z),$ then $(\bar{x},\bar{y})$ is a stationary point of the problem (\ref{111.1}) with a multiplier $(\bar{\lambda}^e, \bar{\alpha}^z, \bar{\gamma}^z)$, namely,
\begin{equation}\label{2.3xxx}
\left\{
\begin{array}{ll}
0\in\nabla f(\bar{x})- A^{\top}\bar{\lambda}^e+C^{\top}(\bar{\gamma}^z-\bar{\alpha}^z)+ N_{\mathcal{X}}(\bar{x}),\\
0\in \nabla \theta(\bar{y})- F^{\top}\bar{\lambda}^e+D^{\top}(\bar{\gamma}^z-\bar{\alpha}^z)+ N_{\mathcal{Y}}(\bar{{y}}),\\
0\leq \bar{\alpha}^z \perp (C\bar{x}+D\bar{y}-r)\geq 0,\ 0\leq \bar{\gamma}^z \perp (s-(C\bar{x}+D\bar{y}))\geq 0,\\
A\bar{x}+B\bar{y}=b.
\end{array}
\right.
\end{equation}
}
\end{lemma}

For the current iteration point $(x_{k},y_{k},z_{k})$ satisfying $x_{k}\in \mathcal{X},~y_{k}\in \mathcal{Y}$ and $r\leq~z_{k}\leq s$,
based on the splitting subproblems (\ref{1.6}) and (\ref{1.7}), we consider the following two subproblems:

%
%we consider the QP subproblem of (\ref{111.2})
%\be\label{111.3}
%\ba{ll} \min & \nabla f(x_{k})^{\top}(x-x_{k})+\frac{1}{2}\|x-x_{k}\|^{2}_{H_{k}^{x}}+\nabla\theta(y_{k})^{\top}(y-y_{k})+\frac{1}{2}\|y-y^{k}\|^{2}_{H_{k}^{y}}\\
%{\rm s.t.}&Ex+Fy+Gz=c,\\
%&x\in \mathcal{X},~y\in \mathcal{Y},~z\geq0.
%\ea
%\ee
%
%The ADMM splitting idea is applied to the problem (\ref{111.3}) to form the following three subproblems:
\be\label{111.4}
 \min\limits_{x\in\mathcal{X}}\ \nabla f(x_{k})^{\top}(x-x_{k})+\frac{1}{2}\|x-x_{k}\|^{2}_{H_{k}^{x}}+\frac{\beta}{2}\|Ex+Fy_{k}+Gz_{k}-c-\frac{\lambda_{k}}{\beta}\|^{2}
\ee
and
\be\label{111.5}
 \min\limits_{y\in\mathcal{Y}}\ \nabla \theta(y_{k})^{\top}(y-y_{k})+\frac{1}{2}\|y-y_{k}\|^{2}_{H_{k}^{y}}+\frac{\beta}{2}\|Ex_{k}+Fy+Gz_{k}-c-\frac{\lambda_{k}}{\beta}\|^{2}.
\ee
%
%\be\label{111.6} \ba{ll} \min\ \frac{\beta}{2}\|Ex_{k}+Fy_{k}+G-c-\frac{\lambda_{k}}{\beta}\|^{2}\\\
%{\rm s.t.~}\  z\geq 0. \ea
%\Leftrightarrow \ba{ll} \min\ \frac{\beta}{2}\|z-\hat{z}_{k+1}\|^{2}\\\
%{\rm s.t.~}\  z\geq 0,
%\ea
%\ee
%where the vector $\hat{z}_{k+1}$ is defined by (\ref{x3}).

In view of $\mathcal{X}$ and $\mathcal{Y}$ being nonempty closed convex sets, by \cite[Corollary 3.4.2]{jianjinbao_2010}, we know that subproblems~(\ref{111.4})~and ~(\ref{111.5}) have unique optimal solutions~$\tilde{x}_{k+1}$ and $~\tilde{y}_{k+1}$ under Assumption \ref{ass2}, respectively.
Furthermore, it follows from \cite[Theorem 6.12]{Rockafellar_Wets_1998} that the optimality conditions of the subproblems above are:
\begin{equation}\label{113.0}
0\in \nabla f(x_{k})+H_{k}^{x}(\tilde{x}_{k+1}-x_{k})+\beta E^{\top}(E\tilde{x}_{k+1}+Fy_{k}+Gz_{k}-c-\frac{\lambda_{k}}{\beta})+N_{\mathcal{X}}(\tilde{x}_{k+1}),
\end{equation}
and
\begin{equation}\label{113.1}
0\in \nabla \theta(y_{k})+H_{k}^{y}(\tilde{y}_{k+1}-y_{k})+\beta F^{\top}(Ex_{k}+F\tilde{y}_{k+1}+Gz_{k}-c-\frac{\lambda_{k}}{\beta})+N_{\mathcal{Y}}(\tilde{y}_{k+1}),
\end{equation}
%\be\label{113.2}\left\{ \ba{ll}
%\beta G^{\top}(Ex_{k}+Fy_{k}+G\tilde{z}_{k+1}-c-\frac{\lambda_{k}}{\beta})-\alpha^{z}_{k}=0,\\
%0\leq\alpha^{z}_{k}\perp\tilde{z}_{k+1}\geq0.
%\ea\right.
%\Longleftrightarrow\left\{ \ba{ll} \beta(\tilde{z}_{k+1}-\hat{z}_{k+1})-\alpha^{z}_{k}=0\\\
%0\leq\alpha^{z}_{k}\perp\tilde{z}_{k+1}\geq0.
%\ea
%\right.
%\ee
% where~$N_{\mathcal{X}}(\tilde{x}_{k+1})$ and $N_{\mathcal{Y}}(\tilde{y}_{k+1})$ are the normal cones of closed convex sets $\mathcal{X}$ and $\mathcal{Y}$ at the optimal solutions $\tilde{x}_{k+1}$ and $\tilde{y}_{k+1}$, respectively.
In a fashion similar to (\ref{2.4}), we also define the direction $(d_{k}^{x},d_{k}^{y})$ by the optimal solutions $\tilde{x}_{k+1}$ and $~\tilde{y}_{k+1}$ of problems (\ref{111.4})~and~(\ref{111.5}), respectively.
From (\ref{2.4}), (\ref{y1}), (\ref{113.0}) and (\ref{113.1}), one can obtain
\be\label{111.9}
-\nabla_{x}\mathcal{L}_{\beta}(w_{k})-(H_{k}^{x}+\beta E^{\top}E)d_{k}^{x}\in N_{\mathcal{X}}(\tilde{x}_{k+1}),
\ee
\be\label{111.10}
-\nabla_{y}\mathcal{L}_{\beta}(w_{k})-(H_{k}^{y}+\beta F^{\top}F)d_{k}^{y}\in N_{\mathcal{Y}}(\tilde{y}_{k+1}).
\ee
Since $\mathcal{X}$ and $\mathcal{Y}$ both are convex, the optimality condition (\ref{111.9}) and (\ref{111.10}) can be rewritten as % (let~$x=x_{k}\in\mathcal{X}$),
(we refer the interested readers to \cite[Theorem 6.12]{Rockafellar_Wets_1998} for more details)
\be\label{111.9b}
(-\nabla_{x}\mathcal{L}_{\beta}(w_{k})-(H_{k}^{x}+\beta E^{\top}E)d_{k}^{x})^\top (x-\tilde{x}_{k+1})\leq 0,\ \forall\ x\in \mathcal{X},
\ee
\be\label{111.10b}
(-\nabla_{y}\mathcal{L}_{\beta}(w_{k})-(H_{k}^{y}+\beta F^{\top}F)d_{k}^{y})^\top (y-\tilde{y}_{k+1})\leq 0,\ \forall\ y\in \mathcal{Y}.
\ee
Now, letting $x=x_k$ and $y=y_k$ in the inequalities \eqref{111.9b} and \eqref{111.10b}, respectively, and in view of $d^x_k=\tilde{x}_{k+1}-x_k$ and $d^y_k=\tilde{y}_{k+1}-y_k$, we have
\[
\nabla_{x}\mathcal{L}_{\beta}(w_{k})^{\top}d_{k}^{x}\leq
-\|d_{k}^{x}\|^{2}_{(H_{k}^{x}+\beta E^{\top}E)},\ \nabla_{y}\mathcal{L}_{\beta}(w_{k})^{\top}d_{k}^{y}\leq-\|d_{k}^{y}\|^{2}_{(H_{k}^{y}+\beta F^{\top}F)}.
\]
Obviously, this, together with (\ref{3.4}) and (\ref{3.5}), we have
\be\label{1.9}
\nabla_{u}\mathcal{L}_{\beta}(u_{k},\lambda_{k})^{\top}d_{k}^{u}\leq-\|d_{k}^{u}\|^{2}_{H_{k}^{u}}.
\ee
The inequality above shows that $\mathcal{L}_{\beta}(\cdot,\lambda_{k})$ is monotonously decreasing along direction $d_{k}^{u}$ at $u_{k}$.

Based on the analysis above, an extension of the previous the GLC-S-SQP algorithm A is proposed as follows.

{\bf GLC-S-SQP algorithm B}

\indent{\bf{Step 0}}~(Initialization) It is the same as {\bf{Step 0}} of the GLC-S-SQP algorithm A except that the initial iteration point  $x_{0}\in \mathcal{X},~y_{0}\in \mathcal{Y}$.
%Choose parameters~$\rho,\sigma\in(0,1),~\beta,\xi>0$ and initial iteration point $w_{0}:=(u_0,\lambda_{0})=(x_{0},y_{0},z_{0},\lambda_{0}^{e},\lambda_{0}^{ie})$ satisfying:~$x_{0}\in \mathcal{X},~y_{0}\in \mathcal{Y},~r\leq z_{0}\leq s,$ two symmetric matrices $H_{0}^{x}\in \mathfrak {R}^{n_{1}\times n_{1}}$ and $H_{0}^{y}\in \mathfrak {R}^{n_{2}\times n_{2}}$ such that~$H_{0}^{x}+\beta E^{\top}E\succ 0,~H_{0}^{y}+\beta F^{\top}F\succ 0.$ Set~$k=0.$

{\bf{Step 1}}~(Solving subproblems)
Solving the two~${\rm{QP}}$ subproblems~(\ref{111.4}) and (\ref{111.5}) to generate the (unique) optimal solutions ~$\tilde{x}_{k+1}$ and $\tilde{y}_{k+1},$ respectively. And~$\tilde{z}_{k+1}$ is generated by (\ref{x3})--(\ref{3.9}).

{\bf{Step 2, Step 3}} and {\bf Step 4} are similar to the associated steps in the GLC-S-SQP algorithm A.

\begin{remark}{\rm
Whether the algorithm above can be implemented effectively depends on whether the two subproblems (\ref{111.4}) and (\ref{111.5}) can be solved  effectively. In particular, if the two closed convex sets $\mathcal{X}$ and $\mathcal{Y}$ are both affine manifolds, the two subproblems (\ref{111.4}) and (\ref{111.5}) can be reduced as standard QP, then they can be solved efficiently.
 }
\end{remark}

%{\bf GLC-ADMM-SQP Algorithm~{\rm $II_{+}$}(the extension of GLC-ADMM-SQP Algorithm~{\rm II}):}
%
%{\bf{Step 0, Step 1, Step 2}} are similar with GLC-ADMM-SQP Algorithm~{\rm $I_{+}$}. {\bf{Step 3, Step 4}} are similar with~GLC-ADMM-SQP Algorithm~{\rm II}.

For the algorithm above, we give the following convergence result.
\begin{theorem}\rm{
Suppose that Assumption \ref{ass2} holds. %Then one of the two situations is true.
%
%(i) If the GLC-S-SQP algorithm B terminates within a finite number of iterations, then it generates a stationary point of
%the original problem (\ref{111.1}).
If the GLC-S-SQP algorithm B generates an infinite sequence $\{w_{k}\}$ of points, and $w_{*}=(x_*,y_*,z_*,\lambda_*:=({\lambda}_*^{e},{\lambda}_*^{ie}))$
is a accumulation point of $\{w_k\}$, then $(x_*,y_*)$ is a stationary point of the original problem (\ref{111.1}),
% and there exists an infinite subsequence $\{(x_{k}, y_{k}, \lambda_{k}^e,\alpha_{k}^{z},\gamma_{k}^{z})\}_{\mathcal{K}}$ converges the stationary point $(x_{*}, y_{*}, \lambda_{*}^e,\alpha_{*}^{z},\gamma_{*}^{z})$,
i.e., the GLC-S-SQP algorithm B is globally convergent.
%Furthermore, there exists an infinite subsequence $\{(\lambda_{k},\alpha_{k}^{z},\gamma_{k}^{z})\}_{\mathcal{K}}$ of multipliers such that it converges to the KKT multiplier associated with the KKT point $(x_*,y_*)$ of (\ref{111.1}).
}\end{theorem}

{\bf Proof}~ In a similar fashion to Lemma \ref{lem2}, if subsequence $\{w_{k}\}_{\mathcal{K}}$ of $\{w_{k}\}$ is bounded, we obtain that the corresponding subsequence $\{d_{k}^u\}_{\mathcal{K}}$, $\{\tilde{u}_{k+1}\}_{\mathcal{K}}$, $\{w_{k+1}\}_{\mathcal{K}}$ and ~$\{(\alpha_{k}^{z},\gamma_{k}^{z}\}_{\mathcal{K}}$ are also bounded, under the Assumption \ref{ass2}. Therefore, the limit $d_{*}=0$ defined by (\ref{114.1}) is zero, so~(\ref{1.10}) also holds. Note that normal cone mapping is closed, then taking the limit in the optimality conditions (\ref{113.0}), (\ref{113.1}) and (\ref{3.0}) for $k\in\mathcal{K}$, we obtain
$$
\left\{\begin{array}{l}\nonumber
-\nabla f(x_{*})+ E^{\top}\lambda_{*}\in N_{\mathcal{X}}({x}_{*}),\\
-\nabla \theta(y_{*})+ F^{\top}\lambda_{*}\in N_{\mathcal{Y}}({y}_{*}),\\
-G^{\top}\lambda_{*}-\alpha_{*}^{z}+\gamma_{*}^{z}=0,\\
0\leq\alpha_{*}^{z}\perp({z}_{*}-r)\geq 0,~0\leq\gamma_{*}^{z}\perp(s-{z}_{*})\geq 0,\\
Ex_{*}+Fy_{*}+Gz_{*}-c=0.
\end{array} \right.
$$
These show that $(x_{*},y_{*},z_{*})$ with the corresponding multiplier $(\lambda_{*},\alpha_{*}^{z},\gamma_{k}^{z})$ satisfies (\ref{4.2x}). Moreover, by Lemma \ref{lem1xxx}, one knows that $(x_*,y_*)$ is a stationary point of the problem (\ref{111.1}), and the theorem is proved.\hfill$\Box$

\section{Applications}

In this section, the numerical validity of our proposed algorithm is tested
by solving a kind of practical economic dispatch problem of power system. The numerical experimental platform is MATLAB R2016a, Intel (R) Core (TM) i5-8500 CPU 3.00GHz RAM 8 GB, Windows 10 (64bite).
\subsection{Problem description}
The economic dispatch (ED) model is a power dispatch (power generation) scheme
that seeks the minimum total power generation cost of a power supply system under the
physical and system constraints of the unit, and under the status that the start and stop
states of the unit set are determined, more details can be found in \cite{WWs,Theerthamalai A_Maheswarapu S_2010}. Its mathematical model can be described as follows:

(1) The objective function of ED is
\begin{equation}\label{5.1}
\min~ F_{c}(p)=\sum\limits_{i=1}^{N}\sum\limits_{t=1}^{\top}(a_{i}p_{_{i,t}}^{3}+b_{i}p_{_{i,t}}^{2}+c_{i}p_{_{i,t}}+d_{i}),
\end{equation}
where~$p_{_{i,t}}$ is the output variable of the unit $i$ in the period $t$, $a_{i},~b_{i},~c_{i},~d_{i}$ are the cost function coefficients of unit $i$, $T$ is the number of optimization periods and $N$ is the number of units.

(2) The constraint conditions  can be defined as follows.

   \indent The power balance constraint:
 \begin{equation}\label{5.2}
 \sum\limits_{i=1}^{N}p_{_{i,t}}=p_{_{D,t}},~t\in\{1,2,\ldots,T\},
 \end{equation}
where~$P_{_{D,t}}$ is the whole network load of period $t$.

\indent The upper and lower output constraint:
\begin{equation}\label{5.3}
p_{_{i,\min}}\leq p_{_{i,t}}\leq p_{_{i,\max}},~i\in\{1,2,\ldots,N\},~t\in\{1,2,\ldots,T\},
\end{equation}
where~$P_{i,\min}~(P_{i,\max})$ is the minimum (maximum) output for unit $i$.

 \indent The unit climbing rate constraint:
 \begin{equation}\label{5.4}
 -D_{i}\leq p_{_{i,t}}-p_{_{i,t-1}}\leq U_{i},~i\in\{1,2,\ldots,N\},~t\in\{1,2,\ldots,T\},
 \end{equation}
where~$D_{i}$ and $U_{i}$ are the upper and lower climbing rate constraints of unit $i$. $p_{_{i,0}}$ is the initial power of unit $i$.
For the above inequality \eqref{5.4}, we transform the unit climbing rate constraint from inequality into equality by introducing the slack variable $q_{i, t}$, then (\ref{5.1})-(\ref{5.4}) can be summarized as the following optimization problem:
\be\label{5.5}
\ba{ll} \min &F_{c}(p)=\sum\limits_{i=1}^{N}\sum\limits_{t=1}^{\top}(a_{i}p_{_{i,t}}^{3}+b_{i}p_{_{i,t}}^{2}+c_{i}p_{_{i,t}}+d_{i})\\\
\ {\rm s.t.}&\sum\limits_{i=1}^{N}p_{_{i,t}}=p_{_{D,t}},~~~~~~t\in\{1,2,\ldots,T\},\\
&-p_{_{i,t}}+p_{_{i,t-1}}+q_{_{i,t}}=0,~~~i\in\{1,2,\ldots,N\},~~t\in\{1,2,\ldots,T\},\\
&p_{_{i,\min}}\leq p_{_{i,t}}\leq p_{_{i,\max}},~~~i\in\{1,2,\ldots,N\},~~t\in\{1,2,\ldots,T\},\\
&-D_{i}\leq q_{_{i,t}}\leq U_{i},~~~i\in\{1,2,\ldots,N\},~~t\in\{1,2,\ldots,T\}.
\ea
\ee
The scale of  the ED model (\ref{5.5}) is as follows: the numbers of variables, equality constraints and box constraints are  $n:=2NT$,  $m:=(N+1)T$  and $2NT$, respectively. The scale of (\ref{5.5})  increases rapidly as  $N$ increases. For example, if $N=200$ and $T=24$, the scale $(n,m)=(9600,4824)$.

In order to solve the above problem by using the GLC-S-SQP algorithm A, and considering the characteristics of engineering in the power system
economic dispatch problem, we divide $p_{_{i,t}}$ equally into two parts. Taking~$N_{1}=[\frac{N}{2}],~N_{2}=N-N_{1},$
$p_{1}=(p_{_{i,t}},~i=1,\ldots,N_{1},~t=1,\ldots,T)\in \mathfrak{R}^{N_{1}T},~p_{2}=(p_{_{i,t}},~i=N_{1}+1,\ldots,N,~t=1,\ldots,T)\in \mathfrak{R}^{N_{2}T},$~$q=(q_{i,t})\in \mathfrak{R}^{NT}.$
Thus, the problem (\ref{5.5}) is equivalent to the following form:
\be\label{5.6}
\ba{ll} \min &F_{c}(p_{1},p_{2})=f(p_{1})+\theta(p_{2})\\\
\ {\rm s.t.}&\left(\ba{c} E_{1}\\M_{1}\\ 0_{(N_{2}T\times N_{1}T)} \ea\right)p_{1}+\left(\ba{c} E_{2}\\0_{(N_{1}T\times N_{2}T)}\\M_{2} \ea\right)p_{2}+\left(\ba{c} 0_{(T\times NT)}\\F_{1}\\F_{2} \ea\right)q=\left(\ba{c} p_{_{D}}\\-\hat{q}^{1}_{_{0}}\\-\hat{q}^{2}_{_{0}} \ea\right),\\
&p^{1}_{\min}\leq p_{1}\leq p^{1}_{\max},\\
&p^{2}_{\min}\leq p_{2}\leq p^{2}_{\max},\\
&-D\leq q\leq U,
\ea
\ee
where
$$f(p_{1})=\hat{p}_{1}^{\top}A_{1}{p}_{1}+{p}_{1}^{\top}B_{1}{p}_{1}+c{1}{p}_{1}+d{1},~\theta(p_{2})=\hat{p}_{2}^{\top}A_{2}{p}_{2}+{p}_{2}^{\top}B_{2}{p}_{2}+c{2}{p}_{2}+d{2},$$
$$\hat{p}_{1}=(p^{2}_{_{i,t}},i=1,\ldots,N_{1},t=1,\ldots,T)\in \mathfrak{R}^{N_{1}T},$$$$\hat{p}_{2}=(p^{2}_{_{i,t}},i=N_{1}+1,\ldots,N,t=1,\ldots,T)\in \mathfrak{R}^{N_{2}T},
E_{1}=(I_{\top},\ldots,I_{\top})\in\mathfrak{R}^{T\times N_{1}T},$$
$$~E_{2}=(I_{\top},\ldots,I_{\top})\in\mathfrak{R}^{T\times N_{2}T},\ F_{1}=(I_{N_{1}T},0_{N_{1}T\times N_{2}T}),$$
$$F_{2}=(0_{N_{2}T\times N_{1}T},I_{N_{2}T}),\hat{q}^{1}_{_{0}}=(p_{_{1,0}},0,\ldots,0,p_{_{2,0}},0,\ldots,0,p_{_{N_{1},0}},0,\ldots,0)^{\top},$$
$$\hat{q}^{2}_{_{0}}=(p_{_{N_{1}+1,0}},0,\ldots,0,p_{_{N_{1}+2,0}},0,\ldots,0,p_{_{N,0}},0,\ldots,0)^{\top},\ p_{_{D}}=(p_{_{D,1}},\ldots,p_{_{D,T}})^{\top},$$
$$A_{1}={\rm diag}(a_1I_T,a_2I_T,\ldots,a_{N_{1}}I_T),~A_{2}={\rm diag}(a_{N_{1}+1}I_T,\ldots, a_{N}I_T)$$$$\ B_{1}={\rm diag}(b_1I_T,\ldots, b_{N_{1}}I_T), B_{2}={\rm diag}(b_{N_{1}+1}I_T,\ldots, b_{N}I_T),$$$$ c{1}=(c_{1}e_{\top},\ldots,c_{N_{1}}e_{\top}),~c{2}=(c_{N_{1}+1}e_{\top},\ldots,c_{N_{1}}e_{\top}),$$
$$d1=T\sum\limits_{i=1}^{N_{1}}d_{i},\ d2=T\sum\limits_{i=N_{1}+1}^{N}d_{i},\ D=(D_{1},\ldots,D_{1},\ldots,D_{N},\ldots,D_{N})^{\top}\in \mathfrak{R}^{NT},$$$$U=(U_{1},\ldots,U_{1},\ldots,U_{N},\ldots,U_{N})^{\top}\in \mathfrak{R}^{NT},$$
$$M_{1}={\rm diag}(M_{0},M_{0},\ldots,M_{0})\in\mathfrak{R}^{N_{1}T\times N_{1}T},$$$$M_{2}={\rm diag}(M_{0},M_{0},\ldots,M_{0})\in\mathfrak{R}^{N_{2}T\times N_{2}T},$$
$$M_{0}=\left(
 \begin{array}{llll}
 -1& \ & \ & \\
 1 & -1 & \ &  \\
  \ & \ddots & \ddots &\\
  & &   1 & -1
 \end{array}
 \right)_{T\times T},$$
and $e_{\top}=(1,1,\ldots,1)$~is the $T$ dimension row vector, $I_T,~I_{N_{1}T}$ and $I_{N_{2}T}$~are $T$,~$~N_{1}T$ and $N_{2}T$ order  identity matrices, respectively.

\subsection{Numerical results and analysis}

In this subsection, by copying the data of 5 units, we generate 20
ED instances, and their structures are shown in Table 1, we refer the interested readers to \cite{Theerthamalai A_Maheswarapu S_2010} for more details.
 By solving this subclass of
the ED instances, we compare the GLC-S-SQP algorithm A with the famous OPTI
solver with version~2.28 downloaded from https://github.com/jonathancurrie/OPTI/releases and an augmented-Lagrange-based SQP
algorithm for the general linear constrained two-block nonconvex optimization problem (\ref{1.2}) (GLC-AL-SQP for short).
Now, we briefly describe the steps of the GLC-AL-SQP algorithm for the problem (\ref{1.2}) as follows.

\newpage

%we solve these instances by GLC-ADMM-SQP, OPTI and GLC-AL-SQP. The numerical results are shown in Table 5-2, where $F_{c}(P^{\ast})$ represents the obtained optimal value, C$_{\rm t}$ represents the CPU calculation time (seconds), and
%RE represents the relative error (RE) of the optimal values with OPTI; for example,
%\[
%{\rm RE}=\frac{F_{c}(P^{\ast})_{({\rm obtained \ by\ GLC-ADMM-SQP})}-F_{c}(P^{\ast})_{({\rm obtained \ by\ OPTI})}}{F_{c}(P^{\ast})_{({\rm obtained \ by\ OPTI})}}\times 100\%.
%\]

\begin{center} {{ Table~1:} The structures of 20 instances obtained by copying the
5-unit system}
\end{center}\setlength{\tabcolsep}{1.6mm}
\begin{longtable}{c c  c c c c c c c c c  c  c c c c c c   c}
\hline
\cline{1-19}
\hline
\multirow{2}{*}{No.}& &\multicolumn{5}{c}{Unit}& & \multicolumn{1}{c}{\multirow{2}{*}{ N}}& &\multirow{2}{*}{No.}& &\multicolumn{5}{c}{Unit}& &\multicolumn{1}{c}{\multirow{2}{*}{N}}\\
\cline{3-7}\cline{13-17}
\linespread{1}\selectfont
   & &  1  & 2  & 3  & 4  & 5  & &   & &   & &  1  & 2  & 3  & 4  & 5  & &  \\
 \hline\cline{1-19}
 1 & &  1  & 2  & 3  & 2  & 2  & &10 & & 11& &	20 & 24 & 27 & 20 &	19 & &110\\
 2 & &  3  & 3  & 3  & 3  & 3  & &15 & & 12& &	22 & 26 & 29 & 22 &	21 & &120\\
 3 & &  4  & 4  & 4  & 4  & 4  & &20 & & 13& &	26 & 30 & 30 & 22 &	22 & &130\\
 4 & & 	5  & 6  & 7  & 7  & 5  & &30 & & 14& &	30 & 33 & 32 & 25 &	30 & &150\\
 5 & &	5  & 10 & 10 & 5  & 10 & &40 & & 15& &	34 & 37 & 36 & 29 &	34 & &170\\
 6 & &	8  & 11 & 12 & 9  & 10 & &50 & & 16& &	36 & 39 & 38 & 30 &	37 & &180\\
 7 & &	10 & 14 & 16 & 15 & 15 & &70 & & 17& &	40 & 44 & 41 & 34 &	41 & &200\\
 8 & &	13 & 18 & 18 & 13 & 18 & &80 & & 18& &	44 & 48 & 45 & 38 &	45 & &220\\
 9 & &	12 & 20 & 25 & 20 & 13 & &90 & & 19& &	48 & 52 & 48 & 40 &	52 & &240\\
 10& &	18 & 22 & 25 & 18 &	17 & &100& & 20& &	50 & 54 & 50 & 42 &	54 & &250\\
\hline\cline{1-19}
\end{longtable}

%\begin{algorithm}%[H]
%\caption{(AL-SQO)}\label{AlgIII}

\begin{enumerate}
\setlength{\itemindent}{13em}

\item[\textbf{GLC-AL-SQP algorithm}]

\item[\textbf{Step 0, Step 2 and Step 3:}] These three steps
are the same as the GLC-S-SQP algorithm A. And Steps 1 and 4 are as follows.
\setlength{\itemindent}{3em}
\item[\textbf{Step 1}] Solve the QP subproblem
  \begin{equation}\label{5.15}
  \begin{array}{cl} \min & \nabla f(x_{k})^{\top}(x-x_{k})+\frac{1}{2}\|x-x_{k}\|^{2}_{H_{k}^{x}}+\nabla \theta(y_{k})^{\top}(y-y_{k})+\frac{1}{2}\|y-y_{k}\|^{2}_{H_{k}^{y}}\nonumber\\
&+\frac{\beta}{2}\|Ex+Fy+Gz-c-\frac{\lambda^k}{\beta}\|^{2}
\\  {\rm s.t.} &l\leq x\leq u,
p\leq y\leq q,
r\leq z\leq s,
  \end{array}
  \end{equation}
  to obtain a (unique) optimal solution $(\tilde{x}^{k+ 1},\tilde{y}^{k+1},\tilde{z}^{k+1})$.
\setlength{\itemindent}{3em}
\item[\textbf{Step 4}] Generate two new symmetric matrices ~$H_{k+1} ^{x}$ and~$H_{k+1} ^{y}$ are the symmetric approximation matrices of $\nabla^{2}f(x_{k+1})$ and $\nabla^{2}\theta(y_{k+1})$, and such that the matrix $H^u_{k+1}$ defined by \eqref{3.5b} is positive definite.
    %satisfy $diag(H^{x}_{k+1}+\beta E^{\top}E,H^{y}_{k+1}+\beta F^{\top}F,\beta G^{\top}G)\succ 0$.
    Set $k:=k+1$, and return to Step 1.
%\item[\textbf{Step 4}]  Generate
%  new symmetric matrix  $H_{k+1}:={\rm diag}(
%  H^x_{k+1},H^y_{k+1})$ which is a reasonable approximation of ${\rm diag}(\nabla^2f(x^{k+1}),\nabla^2\theta(y^{k+1})$ and satisfies\\
%  $(H_{k+1}+\beta (A\ B)^{{\top}}(A\ B))\succ0.$
%  Set $k:=k+1$, and return to Step 1.
\end{enumerate}
%\end{algorithm}

In the experimenttal process, all values of parameters $a_i, b_i, c_i, d_i$, $P_{D,t}$, $p_{i,\min}, p_{i,\max}, D_i$ and $U_i$ et al. are chosen as in \cite{Theerthamalai A_Maheswarapu S_2010}, and
$T=24$. The parameters  in the GLC-S-SQP and GLC-AL-SQP algorithms are uniformly chosen as:
$$\rho=0.8,~\xi=0.001,~\beta=2000,~\sigma=0.8,~\lambda={\rm ones}(NT+T,1).$$
Initial iteration points of each instance is selected as $$(p_{1}^{0},p_{2}^{0},q^{0})=(p_{\min}^{1},p_{\min}^{2},-D),$$
we adopt a unified terminated criterion for all problems: $\|d_{k}^u\|_{\infty}\leq 0.005.$  We directly select the Hessian~matrices of the corresponding objective functions as the quadratic coefficient matrices in the QP~subproblems, i.e., $H_{k}^{x}=\nabla^{2}f(x_{k}),~H_{k}^{y}=\nabla^{2}\theta(y_{k}).$ Under the background of power system economic dispatch problem,
the uniformly positive definite of these matrices is always satisfied.

%In order to verify the efficiency of the algorithms, we use~GLC-ADMM-SQP, GLC-SQP~and~OPTI~solver to solve the selected 20 numerical instances on the same experimental platform and software running environment.
The numerical results are shown in Table 2,
where iter represents the number of iterations, $F_{c}(P^{\ast})$ represents
approximate optimal
objective value at the final iteration point, $\varphi_{\rm eq}$ represents $\|Ex^*+Fy^*+Gz^*-c\|_{\infty}$, C$_{\rm t}$ represents the CPU calculation time (seconds), and
RE represents the relative error (RE) of the optimal values with OPTI, for example,
\[
{\rm RE}=\frac{F_{c}(P^{\ast})_{({\rm obtained \ by\ GLC-S-SQP})}-F_{c}(P^{\ast})_{({\rm obtained \ by\ OPTI})}}{F_{c}(P^{\ast})_{({\rm obtained \ by\ OPTI})}}\times 100\%.
\]
For simplicity, denote the ``Sum'' of $\varphi_{\rm eq}$ and C$_t$ of OPTI, GLC-S-SQP and GLC-AL-SQP, respectively.
From the numerical reports in Table 2, we have the following
claims:

\begin{sidewaystable}%[!htbp]
\vspace{102.8mm}
\setlength{\tabcolsep}{120.8mm}
%\tiny
\renewcommand\arraystretch{1.2}
\leftline{\small{\bf Table 2} Numerical results of ED instances obtained by
OPTI, GLC-S-SQP and GLC-AL-SQP
}
\label{T5-2}       % Give a unique label
\setlength{\tabcolsep}{2.8mm}
{ \footnotesize
\begin{tabular}{ccc|ccccc|ccccc}
\noalign{\smallskip}\hline
\multirow{2}*{No.}& \multicolumn{2}{c}{OPTI } & \multicolumn{5}{c}{ GLC-S-SQP} & \multicolumn{5}{c}{ GLC-AL-SQP}   \\
\cline{2-13}
%&{Total cost}&\multirow{2}*{C$_{\rm t}$}&{Total cost}&\multirow{2}*{$\varepsilon^{\rm RE-A1}$}&\multirow{2}*{N$_{\rm i}$}&\multirow{2}*{C$_{\rm t}$}& {Total cost}&\multirow{2}*{$\varepsilon^{\rm RE-A2}$}&\multirow{2}*{N$_{\rm i}$}&\multirow{2}*{C$_{\rm t}$}\\
& {$F_{c}(P^{\ast})$}& {C$_{\rm t}$} & iter& {$F_{c}(P^{\ast})$} & {$\varphi_{\rm eq}$} &{C$_{\rm t}$} & {${\rm RE(\%)}$} & iter & {$F_{c}(P^{\ast})$} & {$\varphi_{\rm eq}$}& {C$_{\rm t}$} & {${\rm RE(\%)}$}  \\
\noalign{\smallskip}\hline
1&1243485.20 & 1.53 &28&1243923.71 &0.024&0.55& 0.0353 &   32&1243479.85&0.045&0.93&-0.0004 \\
2&1833617.31 & 2.68 &29& 1834256.70 &0.025&0.76&0.0349& 32& 1833609.97 &0.062&1.62&-0.0004\\
3&2444823.08 & 5.83 &29& 2445520.86 & 0.028&1.06&0.0285 & 32& 2444813.53&0.082& 2.25&-0.0004\\
4&3650316.85 & 28.84 &29& 3651460.87&0.032&1.71&0.0313 &32&3650302.94&0.121&4.31&-0.0004 \\
5&5083735.09 & 29.57 &29& 5085281.80&0.040&2.75&0.0304 &32& 5083715.20&0.173&5.07&-0.0004 \\
6&6192152.16 & 43.50 &30& 6194045.96 &0.047&3.51&0.0306  &32& 6192128.66&0.205&7.29&-0.0004 \\
7&8636967.25 & 464.27 &31& 8639819.89&0.059&5.51&0.0330  &32& 8636934.92&0.283&12.78&-0.0004 \\
8&9973328.82 & 654.49 &31& 9976340.39&0.061&6.89&0.0302 &32&9973291.26&0.329&16.67&-0.0004 \\
9&11035233.41 &1004.45 &31& 11038863.29&0.067&7.76&0.0329 & 32&11035192.23&0.361&18.78&-0.0004\\
10&12291433.08 &1005.87 &31& 12295231.31&0.073&9.44&0.0309  &32&12291387.38&0.401&23.15&-0.0004 	\\
11&13513839.10 &1004.15 &31& 13517999.46&0.080&11.51&0.0308 & 32& 13513788.92&0.440&27.80&-0.0004 	\\
12&14736246.11 &1009.19 &31& 14740765.65&0.091&12.70& 0.0307  &32& 14736191.58&0.479&34.98&-0.0004\\
13&15975567.93 &1009.17 &31& 15980583.45&0.097&15.62&0.0314  &32& 15975508.93&0.518& 42.12&-0.0004 	\\
14&18492204.68 &1009.88 &31& 18497841.64&0.108&19.96&0.0305  &32& 18492136.47&0.599& 53.50&-0.0004  \\
15&20937025.33 &1015.06 &32& 20943355.76&0.117&25.43&0.0302  &32& 20936948.29&0.677& 68.51&-0.0004 	\\
16&22197414.84 &1015.32 &32& 22204141.87&0.123&28.07&0.0303  &32& 22197333.12&0.719&77.31&-0.0004 	\\
17&24659160.06 &1015.29 &32& 24666551.28&0.136&33.96&0.0300  &32& 24659069.45&0.797& 95.16&-0.0004 	\\
18&27103981.13 &1014.70 &32& 27112072.66&0.149&43.24&0.0299 &32& 27103881.69&0.875&122.02&-0.0004  \\
19&29641688.17 &1015.93 &32& 29650513.94&0.162&50.67&0.0298 &32& 29641579.31&0.958&143.98&-0.0004\\
20&30864098.29 &1021.03 &32& 30873273.98&0.169&55.25&0.0297 &32& 30863985.00&0.997&153.96&-0.0004  \\
\noalign{\smallskip}\hline
Sum & --& 13370.75& --& --& 1.688&336.35 & --&-- &--& 9.121& 912.19 & -- \\
\noalign{\smallskip}\hline
\end{tabular}
}
\end{sidewaystable}

(i) The OPTI solver needs more calculation time to solve 20 ED instances, especially when the number of units exceeds 80, the calculation time has exceeded 1000 seconds, which is unreasonable for solving economic dispatching problems in the real situation. From the perspective of saving calculation time cost, the GLC-S-SQP algorithm A has obvious advantages that 20 examples can be solved effectively in 60 seconds.

(ii) Compared with~OPTI, the RE of the~GLC-S-SQP algorithm A are about 3/10000. This result is still acceptable under the premise of fully considering the saving calculation time cost. In addition, as the scale of the problem increases, the calculation time of the GLC-S-SQP algorithm A are relatively stable, which reflects the good robustness of our proposed algorithm.

(iii) It is found that the GLC-AL-SQP algorithm is superior to the GLC-S-SQP algorithm A in terms of the value of
RE. However, the Sum of $\varphi_{\rm eq}$ and C$_t$ for the GLC-AL-SQP algorithm is inferior to the latter.

Therefore, we preliminarily conclude that
the GLC-S-SQP algorithm A is superior to the OPTI solver and the GLC-AL-SQP algorithm in terms of computing time and computing accuracy.

\section{Conclusions}

In this work, based on the ideas of the splitting algorithms and SQP methods, and by
means of Armijo line search technology with an augmented Lagrangian merit function,
we design a monotone splitting SQP algorithm for solving nonconvex two-block
optimization problems with linear equality, inequality and box constraints. We analyze
the global convergence of the proposed algorithm. In addition, the box constraints are extended to general nonempty closed convex sets. The global convergence of the two algorithms has been proved.
By solving the mid-to-large-scale economic dispatch instances in power systems,
the numerical results show that the proposed algorithm is promising.

We think along with the idea of this work, there are still some interesting and
meaningful problems worth further studying and exploring:

(i) Extend the proposed algorithms to multi-block nonconvex optimization problems.

(ii) Study a Peaceman-Rachford splitting SQP algorithm for two-block optimization with the linear equality
and inequality constraints. Further, explore that the slack variable is not yielded by Armijo line search, but updated by an explicit correction.

(iii) Explore a GLC-S-SQP algorithm with the reasonable iteration complexity and superlinear convergence rate.

\section*{Compliance with ethical standards}

{\bf Conflict of interest} The authors declare that they have no conflict of interest.

\section*{Authors' contributions}

Jinbao Jian carried out the idea of this paper and proposed the description of GLC-S-SQP algorithm A. Guodong Ma carried out the extension of GLC-S-SQP algorithm A and analyzed the global convergence of two algorithms. Xiao Xu drafted the manuscript. Daolan Han carried out the numerical experiments. All authors read and approved the final manuscript.

%% BibTeX users please use one of
%%\bibliographystyle{spbasic}      % basic style, author-year citations
%%\bibliographystyle{spphys}       % APS-like style for physics
%\bibliography{template}   % name your BibTeX data base


%\bibliographystyle{bmc-mathphys} % Style BST file (bmc-mathphys, vancouver, spbasic).
%\bibliography{bmc_article}      % Bibliography file (usually '*.bib' )

%
%

\end{document}